\documentclass[12pt]{amsart}
\usepackage{amsmath}
\usepackage{amscd}
\usepackage{graphics}
\usepackage{latexsym}
\usepackage {amssymb, mathrsfs, graphicx, color, url,epsf}

\def\C{{\mathbb C}}

\def\P{{\mathbb P}}

\begin{document}

\title[]{ The Dolgachev Surface \\ $\;$ \\  \small{\it{ D\lowercase{isproving} H\lowercase{arer}-K\lowercase{as}-K\lowercase{irby} C\lowercase{onjecture}}}}

\author{Selman Akbulut }
\thanks{The author is partially supported by NSF grant DMS 9971440}
\address{Department  of Mathematics, Michigan State University,  MI, 48824}
\email{akbulut@math.msu.edu }
\date{\today}
\begin{abstract} 
We prove that the Dolgachev surface $E(1)_{2,3}$ admits a handlebody decomposition without $1$- and $3$- handles, and we draw the explicit picture of this handlebody. We also locate a ``cork'' inside of  $E(1)_{2,3}$, so that  $E(1)_{2,3}$ is obtained from $E(1)$ by twisting along this cork. 
\end{abstract}
\maketitle

\setcounter{section}{-1}

\vspace{-.1in}

\section{Introduction}

It is a curious question whether an exotic copy of a smooth simply connected $4$-manifold admits a handle decomposition without $1$- and $3$- handles?. Clearly if exotic $S^4$ and $\C\P^2$ exist, their handle decomposition must contain either $1$- or $3$-handles.  Hence It is a particularly interesting problem to find smallest exotic manifolds with this property. Twenty two years ago  Harer, Kas and Kirby conjectured that the Dolgachev surface $E(1)_{2,3}$, which is an exotic copy of $\C\P^2\# 9\bar{\C\P}^2$, must contain $1$- and $3$-handles \cite{hkk}. Recently, Yasui constructed an exotic $\C\P^2\# 9\bar \C\P^2$ with the same Seiberg-Witten invariants as $E(1)_{2,3}$ without $1$- and $3$-handles \cite{y}. Here we disprove Harer-Kas-Kirby conjecture by showing that in fact $E(1)_{2,3}$ itself admits an easy to describe handlebody consisting of only $2$-handles (Figure 41).  As a corollary we show that $E(1)_{2,3}$ admits a simple ``Cork'' as in  the case of some elliptic surfaces E(n) and the Yasui's manifold \cite{ay}. We also show that $E(1)_{2,3}$ is obtained from  $E(1)$ by twisting along this cork. Here I would like to thank Yasui for motivating me to look at this problem.

\vspace{.05in}

Our purpose here is twofold, while disproving this conjecture we also want to fix some consistent conventions. During the construction of the handlebody for $E(1)_{2,3}$, we will set a dictionary between several different descriptions of the elliptic surface $E(1)$ and its exotic copies. In the future we hope to be able use this for constructing of other exotic rational surfaces and Stein fillings.

\vspace{.05in}

Finally, we want to emphasize the simple reoccurring theme in this paper, as in many of our previous works in $4$-manifolds, it is a trivial to state but hard and tedious to practice principle: ``If you keep turning handles of a 4-manifold upside down, while isotoping and canceling, you get a better picture of the manifold''. We have found invoking this principle is often the last saving step, when a proof gets hopelessly stuck during long and hard handle slides. For example in \cite{a1} and \cite{a2} we used this technique in a decisive way, also the proof of \cite{g} was based on \cite{ak} where an arduous turning upside down process  had already been performed. In this paper we use this technique twice. We first apply \cite{a3} to describe a handlebody for $E(1)_{2,3}$ and cancel its $1$-handles, to cancel its $3$-handles we turn this handlebody upside down and cancel the corresponding $1$-handles, then finally by turning it upside down once again we obtain a very simple explicit handlebody for $E(1)_{2,3}$ (no we don't get the same thing when we turn a handlebody upside down twice, since during this process we are also simplifying it by handle slides and cancellations).

\section{$E(1)$}

We start with $\C \P^2 \# 9 \bar{\C\P}^2$, which is also know as the elliptic surface $E(1)$.  It is easy to see that Figure 1 describes a handlebody for $E(1)$, where $\{ h, e_1,...,e_9\}$ corresponds to the standard homology generators of $\C \P^2 \# 9 \bar{\C\P}^2$.

 \begin{figure}[ht]  \begin{center}  
\includegraphics[width=.8\textwidth]{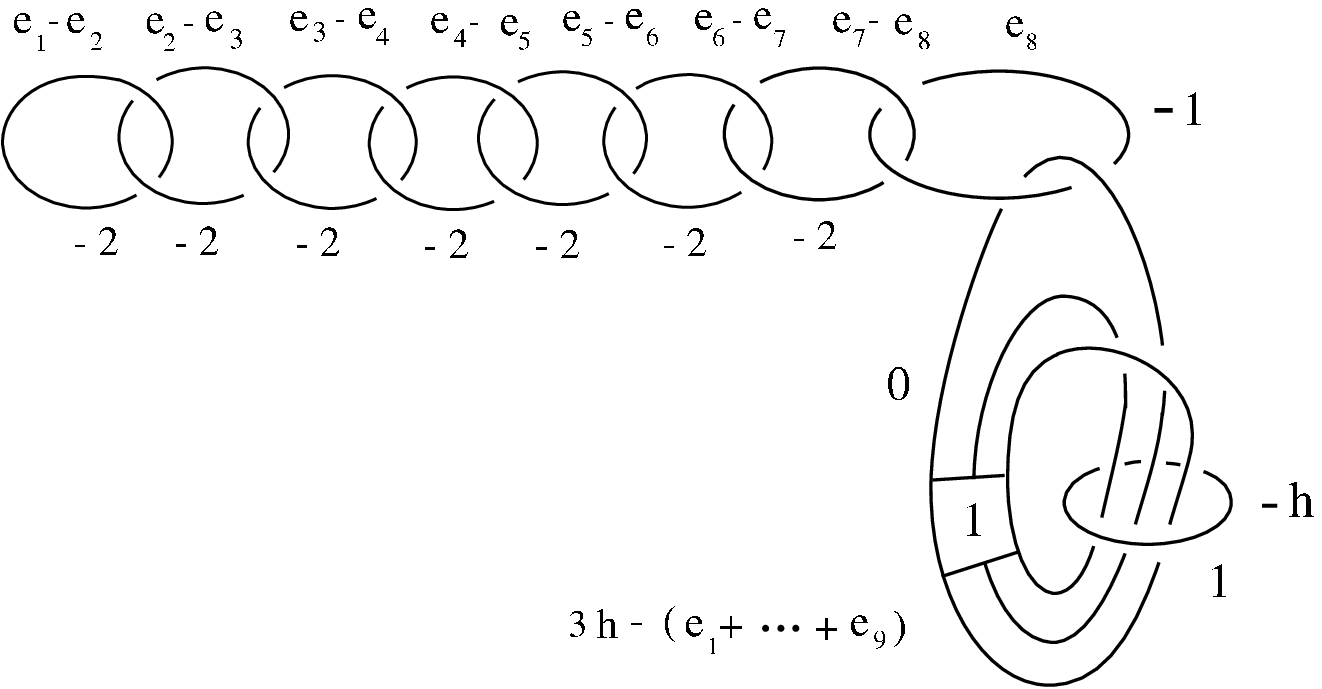}   \caption{}    \end{center}
  \end{figure}

The complex surface $E(1)$ admits a  Lefschetz fibration over $S^2$ with regular torus fibers and 12 singular fishtail fibers, with monodromy $(ab)^6=1$,  where  $a,  b$ are the Dehn twists along  the two standard generators in the mapping class group of $T^2$. The corresponding handlebody of this is given by the first picture in Figure 2.  Alternatively we can use the word $ a^2b^2a^2ba^4b=1$ to describe the same fibration, which is the second picture of Figure 2.  In the both pictures of Figure 2 the unmarked $2$-handles (whose framings are not specified) are attached by $-1$ framings. The homology class $T=3h -(e_{1}+e_{2}+...+e_{9}) $ corresponds to the torus fiber of $E(1)$. Now in the rest of this section we will identify this description with the handlebody of Figure 1, in particular the cusp in Figure 1 will correspond to the singular fiber $T$. An expert reader, who is comfortable with this identification, may skip rest of this section.

\vspace{.05 in}

Starting from Figure 2, via the obvious handle slides (often indicated by arrows in the figures), we obtain Figures 3, 4, ...  and  finally the  first picture of Figure 9. Reader should think of the sequences of these diffeomorphisms (as well as the subsequent similar ones in the paper) as a short movie.

\vspace{.05 in}

We claim that the first picture of Figure 9 is diffeomorphic to $\C \P^2 \# 9 \bar{\C\P}^2$, and the second picture of Figure 9 describes the handles  in terms of the standard homology generators of $\C \P^2 \# 9 \bar{\C\P}^2$. The sequence of handle slides from Figure 9 through  the first picture of the Figure 12 proves this claim;  they describe the precise diffeomorphism to $\C \P^2 \# 9 \bar{\C\P}^2$. By tracing the standard handles of $\C \P^2 \# 9 \bar{\C\P}^2$ backwards to Figure 9,  we see the identification of the second picture of Figure 9. The steps are self explanatory from the figures (the two unknotted zero-framed $2$-handles at the end are cancelled by two $3$-handles).  

\vspace{.05 in}

Having checked the claim,  we return to Figure 9. Now we perform a few obvious handle slides to Figure 9 and obtain the second picture of Figure 12 . Then by doing more handle slides (as indicated by the arrows) we arrive to the pictures in Figure 13, and then obtain  the first picture of Figure 14.  Amazingly,  it is easy to check that the $-2$ framed circle (which is indicated by the arrow  in the picture) is just the unknotted circle with $0$- framing on the boundary of the rest of the handlebody!. Hence we can simply erase this handle from the picture (this corresponds to canceling a $2$ and $3$ handle pair),  and arrive to the second picture of Figure 14, which is a pretty handlebody picture of $E(1)$. It is very easy to see that this picture is diffeomorphic to $\C \P^2 \# 9 \bar{\C\P}^2$ (blow down $+1$ framed handle, and nine $-1$ framed handles, consecutively). The homology generators are indicated in the picture. In fact, we can go one more step by canceling the last  $3$-handle from this handlebody, and  obtain even a simpler picture of $E(1)$ given in the Figure 1 of the introduction, which consisting of just ten $2$-handles, still shoving the cusp inside. 

\vspace{.05 in}

The  first picture of Figure 15 is just a redraw  of  Figure 1 in a more convenient way, and from this by doing the handle slides as indicated by the arrow we obtain the second picture of Figure 15. Then again by  a few more handle slides, as indicated by the arrow, we obtain the first picture of  Figure 16. An isotopy gives the second picture of Figure 16  describing $E(1)$. This picture has a nice feature of not having any $1$-and $3$-handles, and containing the $E_{8}$ plumbing inside.

\section{ $E(1)_{2,3}$}

$E(1)_{p,q}$ is the complex surface obtained by performing two logorithmic transforms of order  $p,q$ (relatively prime integers) on two distinct regular torus fibers of $E(1)$ (c.f. \cite{hkk}).  This operation was  introduced by Kodaira, and studied by Dolgachev. These surfaces have since come to be known as ``Dolgachev Surfaces'', and it was shown by Donaldson that they are fake copies of $E(1)$. Also it is know that $E(1)_{2,3}$ can be identified by the ``{\it Knot surgered }''  copy $E(1)_{K}$ of $E(1)$, by using the trefoil knot $K$ \cite{fs}, \cite{p}. 

\vspace{.05in}

Recall that $E(1)_{K}$ is the manifold obtained from $E(1)$ by replacing a tubular neighborhood $T^2 \times D^2$ of a regular fiber by $(S^{3}-N(K))\times S^1$, where $N(K)$ is the tubular neighborhood of the knot $K$ in $S^3$. Also recall that, in \cite{a3} and \cite{a4} an algorithm of drawing the handlebody picture of a knot surgered manifold, from the handles of the original manifold, was given. By applying this process to the handlebody of $E(1)$ in Figure 16 (with $K$ trefoil knot), we get the first picture of the Figure 17 describing the handlebody of  the Doglachev surface $E(K)_{K}=E(1)_{2,3}$.  Here recall that in a framed link picture, a slice knot with a dot means that the obvious slice disk is removed from the zero handle $B^4$ (this is usually called a {\it slice $1$-handle}); this is just a generalization of unknot with a dot notation, which had been introduced in \cite{a5}.

\section{Canceling $1$-handles of  $E(1)_{2,3}$}

 After converting the slice $1$-handle to a pair of a $1$-handles and a $2$-handle, and by the obvious $2$-handle slides we get the second picture of Figure 17. Clearly in this picture of $E(1)_{2,3}$ all the $1$-handles are cancelled! (the trivially linking $-1$ framed circles, to the $1$-handles, cancels those $1$-handles).

\vspace{.05in}

Now we have no $1$-handles, but a single $3$-handle. This is because Figure 17 has eleven uncanceled $2$-handles, while E(1) has ten homology generators. So the boundary of the Figure 17 is $S^1\times S^2$, which is capped with a $3$ and  $4$ handle pair (i.e. with $S^1\times B^3$). As usual in the figures we don't (and don't need to) draw $3$ and $4$ handles.

\section{Turning $E(1)_{2,3}$ upside down}

Having cancelled $1$-handles of $E(1)_{2,3}$, we will now cancel its $3$-handle. To cancel the $3$-handle we will turn the handlebody of $E(1)_{2,3}$ upside down and cancel the resulting dual $1$-handle. This is done by finding a diffeomorphism from the boundary of the Figure 17 to $\partial (S^1\times B^3)$ and attach the dual $2$-handles (trivially linking zero-framed circles to the un-cancelled $2$-handles of Figure 17) to $S^1\times B^3$ via this diffeomorphism. 

\vspace{.05in}

Now if we draw the dual $2$-handles in the second picture of Figure 17, and carry them to the more convenient first picture of Figure 17, we get the first picture of Figure 18. The dual $2$-handles are drawn in blue. Next we have to find any diffeomorphism from the boundary of Figure 17 to $S^1\times S^2$. During this process no (black) handles can slide over the dual (blue) $2$-handles, but dual $2$-handles can slide over all other handles, and they can slide over each other. In short, we are allowed to change the interior of Figure 17 (black handles) any way we want until it becomes $S^1 \times B^3$, and during this process we cary along the dual (blue handles) along for the ride. 

\vspace{.05in}

In order not to clutter our figures, we won't always specify the framings in the picures when they are obvious from the context, for example the unmarked dual (blue) circles in Figure 18 have all zero framings.

\vspace{.05in}

By turning slice 1-handle to a zero framed knot (i.e.  by replacing the dot with zero framing), then blowing up an unknot, and sliding it over one of the zero framed knots (the ones going through the bottom $1$-handle), and then blowing it down again we obtain the second picture of Figure 18. Then by isotopies and the indicated handle slides we arrive to Figure 19. By blowing up an unknot, sliding over middle zero handle and blowing it down again (several times), and by an isotopy we arrive to Figure 20. Then by the indicated handle slides, and by blowing down a chain of $-1$ framed knots (this changes the $0$ framings on their dual handles to $+1$'s) and by an isotopy we arrive to Figure 21. Now we blow down the top (large) $+1$ framed knot, which turns $+8$ framed trefoil knot to an unknot with $-1$ framing. We then blow down this resulting $-1$ framed unknot; in order not to clutter the picture we indicate this blowing down operation as a canceling handle pair (i.e. the $-1$ framed unknot becomes a $1$-handle with a dual linking $+1$ framed circle which cancels it). This gives us Figure 22. 

\vspace{.05in}

Figure 22 is the picture of $E(1)_{2,3}$ turned upside down. Note that if we ignore the dual (blue) handles, the Figure 22 is just $S^1\times B^3$. To show this in the figure, we circled (black) framed knots which are canceling the three of the four $1$-handles, leaving just a single  $1$-handle which is just $S^1\times B^3$. 

\section{Cancelling $1$-handles of the upside down $E(1)_{2,3}$}

Now we claim that all the $1$-handles Figure 22 can be cancelled. At first glance this is not evident from this picture. To see this, we perform the indicated handle slides and isotopies which takes us from Figure 22 through Figure 25. Finally by sliding the bottom-right $-1$ framed handle two of the $+1$ framed handles (the two small ones that link the middle $1$-handle trivially) we arrive to Figure 26. Now notice that in the Figure 26 every $1$-handle is cancelled by a $2$-handle (which indicated in the figure by encircling their framings). Hence Figure 26 (or equivalently Figure 25) describes a the upside down handlebody of $E(1)_{2,3}$ without $1$ and $3$ handles. Unfortunately, this is not such a pleasant looking handlebody, it is a bit complicated. To improve its image we will turn it upside down one more time in the next section.

\section{  Turning $E(1)_{2,3}$ upside down second time to improve its image}

 To turn Figure 26 upside down, as before we first need to find a diffeomorphism from the boundary of this figure (equivalently from the boundary of Figure 25) to $S^3$:  Figure 27 is obtained from Figure 25 by first canceling the botom-left $1$-handle (with the $-1$ framed $2$-handle), and by turning the top $1$-handle to a zero framed unknot, and then by applying an isotopy. Now by blowing down the two $-1$ circles in the middle, and by blowing down two chains of $+1$ circles we arrive to the first picture of Figure 28. Then by indicated handle slides, isotopies and blowing downs we arrive Figures 29, 28.. and finally to the last picture of Figure 31 which is $S^3$.
 
 \vspace{.05in}

 Now that we found a diffeomorphism to $S^3$ in the last paragraph, we are ready to turn the Figure 26 upside down: In the Figure 32 we drew the dual $2$-handles (indicated by the small red colored circles) of the Figure 26 (we only need to take the duals of the un-cancelled $2$-handles). Then apply this diffeomorphism to $S^3$,  and attach the images of the dual $2$-handles to $B^4$, i.e. we attach $2$-handles to $B^4$ along the images of the red curves. The steps 
 Figure 32 $\mapsto $ Figure 38, is the same as the steps Figure 26 $\mapsto $ Figure 31, except in this case we are carrying the dual $2$-handle circles (red circles) along.

  \vspace{.05in}
  
  Again keep in mind that the dual (red)  $2$-handles can slide over all other handles, and and they can slide over each other, whereas the other handles can not slide over the dual (red) handles. For example, this explains in Figure 36 how the blue handle with a small the red linking circle moved pass red handles from bottom to the top. So the Figure 38 is the final simple picture of $E(1)_{2,3}$ without $1$- and $3$- handles (as indicated in the figure, the lone $1$-handle is canceled by the $-1$ framed $2$-handle). 
  
  There is even a simpler picture of $E(1)_{2,3}$ given in Figure 41 (either pictures) which will be explained in the next section.

 \vspace{.05in}
 
\noindent {\it Remark: A simple corollary of our proof is: $E(1)$ has infinitely many distinct smooth structures without $1$-handles. This is because, the only thing the move in Figure 17 uses is that `the knot $K\#(-K)$ bounds a ribbon disk with two minima' (i.e. a single ribbon move turns $K\#(-K)$ into two unknotted circles). Clearly there are infinitely many such knots $K$, with distinct Alexander polynomials, so $E(1)_{K}$ give examples of distinct smooth manifolds without $1$-handles.}

\section{A cork decomposition of E$(1)_{2,3}$}

Let $W$ be the contractible Stein manifold described in the following figure
   \begin{figure}[ht!]
  \begin{center}   
\includegraphics{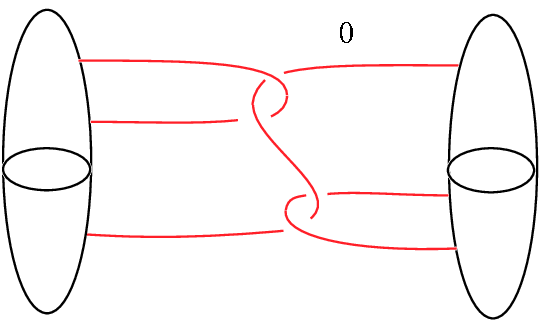}   

{W}    
 \end{center}
   \end{figure}

\noindent and let $f: \partial W\to \partial W$ be the obvious involution, defined by exchanging the positions of a $1$-handle $S^1 \times B^3$  with a $2$-handle $B^2 \times S^2$ in the interior of $W$ (i.e. replacing the ``dot'' and the ``zero'' in the symmetric link of the second picture of Figure 40, which is just an alternative description of $W$). Note that, $W$ is the so called ``positron'' of \cite{am}, and  $\bar{W_{1}}$ in \cite{ay}. Here we claim that $(W, f)$ is a cork of $E(1)_{2,3}$. That is, there is an imbedding $W\subset E(1)_{2,3}$, such that cutting $W$ out and reglueing with the involution $f$ changes the smooth structure of $E(1)_{2,3}$. Let us write $N\cup_{id}W=E(1)_{2,3}$ where $N$ is the complement of $W$. Recall also that $E(1)_{2,3}$ is an irreducible manifold. We will prove our claim by showing a splitting $N\cup_{f}W=P\;\# 5\bar{\C\P}^2$, where $P$ is some smooth $4$-manifold. Later on we will show more advanced version of this namely: $N\cup_{f}W=E(1)$.

\vspace{.05in}

By standard handle slides, from Figure 38 we obtain the two equivalent{\linebreak} diagrams of Figure 39. Then by ignoring some handles, and by the indicated handle slide, we arrive to the second picture of Figure 40 which is $W$. Furthermore, notice that, in the first handlebody diagram of Figure 40, if we replace the``dot'' with zero in the symmetric link (i.e. reglue $W$ with the involution $f$) we get a splitting of  $5$ copies of $\bar{\C\P}^2$. A closer inspection shows that in fact a stronger version of this result holds:

\vspace{.1in}

\noindent {\it {\bf Theorem}:  $E(1)_{2,3}$ is obtained by cork twisting $E(1)$ along the cork $W$.  That is we can decompose   $E(1)_{2,3}=N\cup_{id}W$, so that ${\C\P}^2\# 9\bar{\C\P}^2= N\cup_{f}W$.

\proof First notice that there is even a simpler handlebody of $E(1)_{2,3}$ given in Figure 41 then the one described in Figure 39. To do this just observe that, on the boundary of the first handlebody of Figure 39, the $-2$ -framed circle is isotopic to the $-1$ framed circle of the first handlebody of Figure 41 (which is indicated by the dotted arrow). By a handle slide we obtain the second picture of Figure 41, where $W$ is clearly visible (the circle with dot and the large zero framed circle). Hence the cork twisting of $W$ is given by the first picture of Figure 42 (exchanging ``dot'' with zero framing). Then by the indicated handle slides and diffeomorphism of  Figures 42 to 44 we end up with ${\C\P}^2\# 9\bar{\C\P}^2$. \qed

\newpage

\includegraphics[width=.6\textwidth]{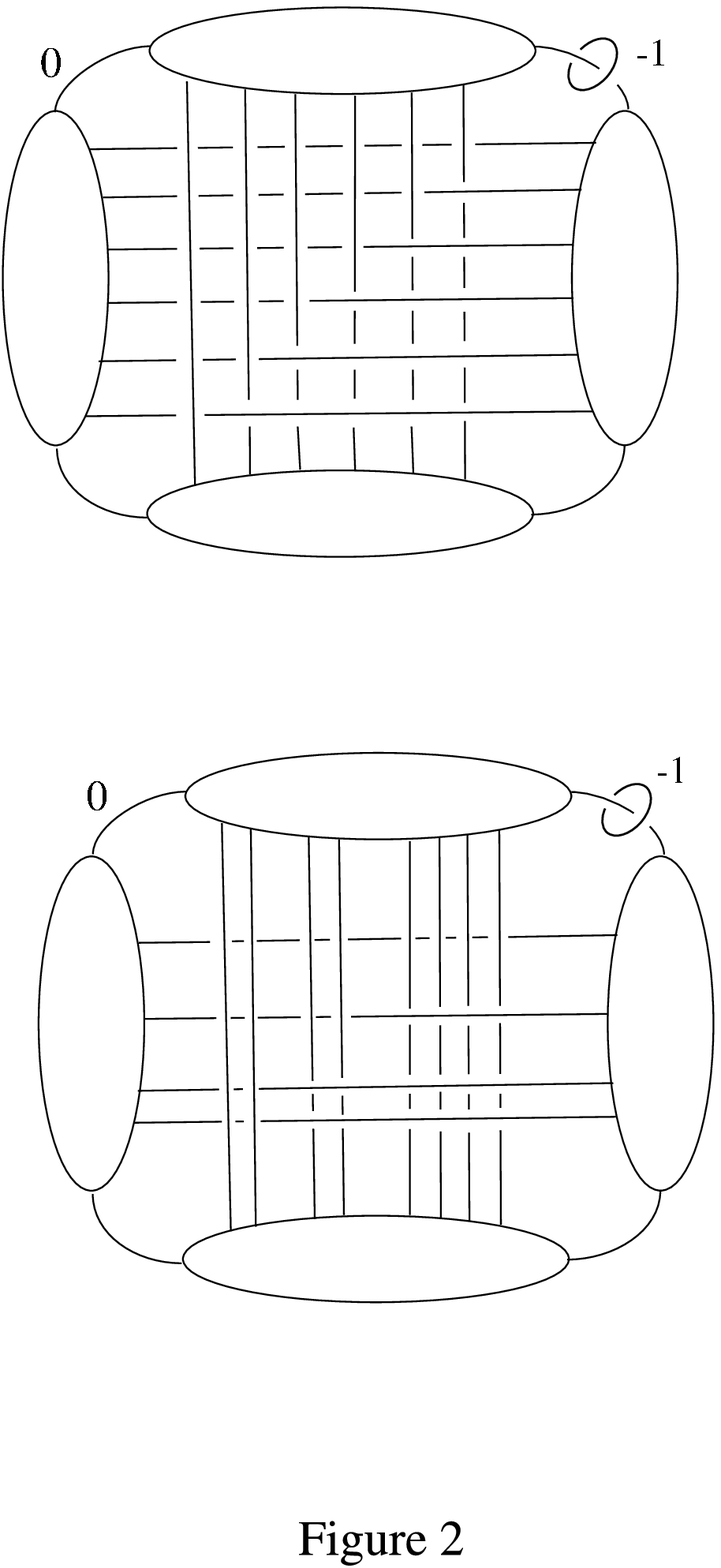}
   
\includegraphics[width=.8\textwidth]{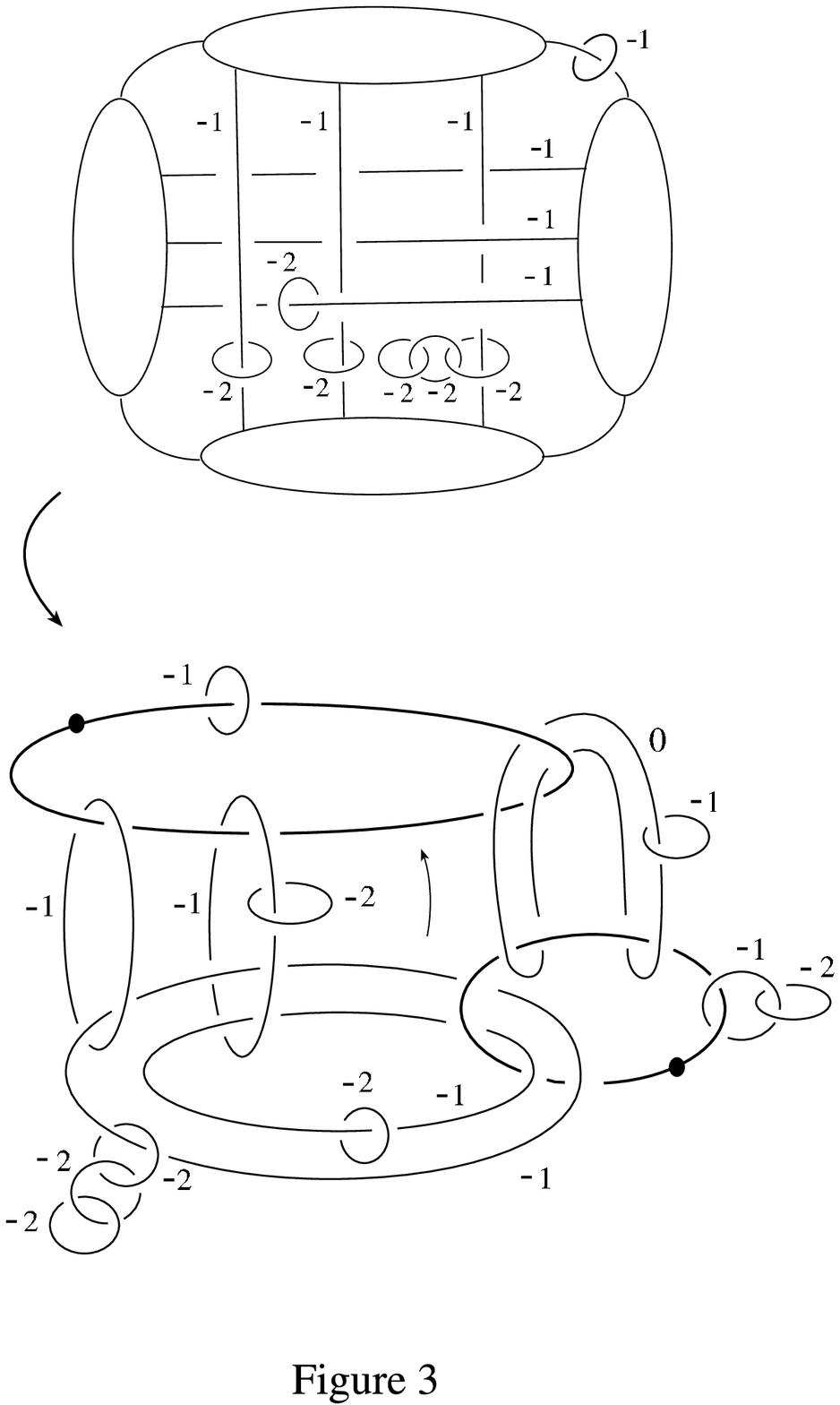}

\includegraphics[width=.9\textwidth]{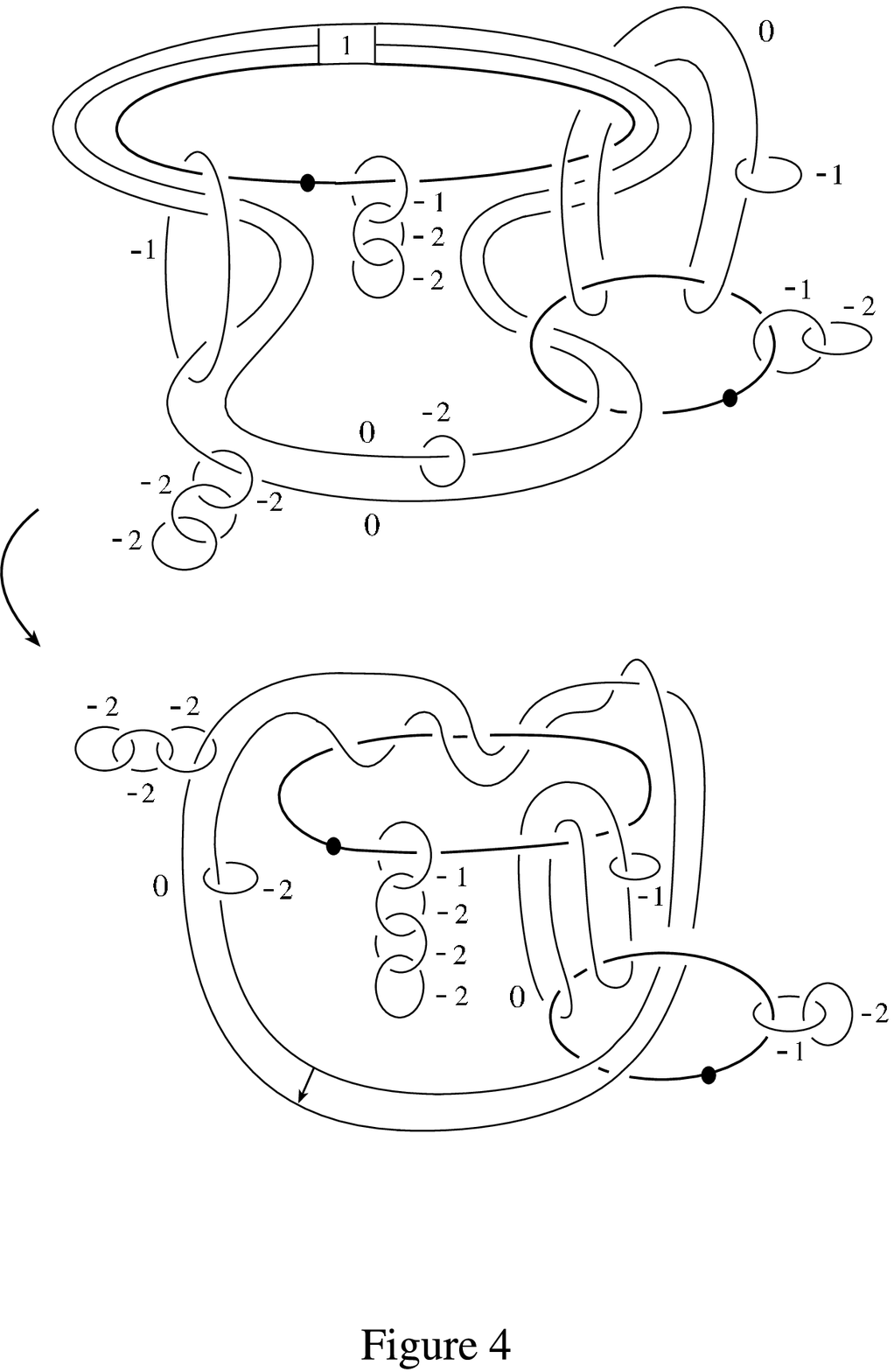}

\includegraphics{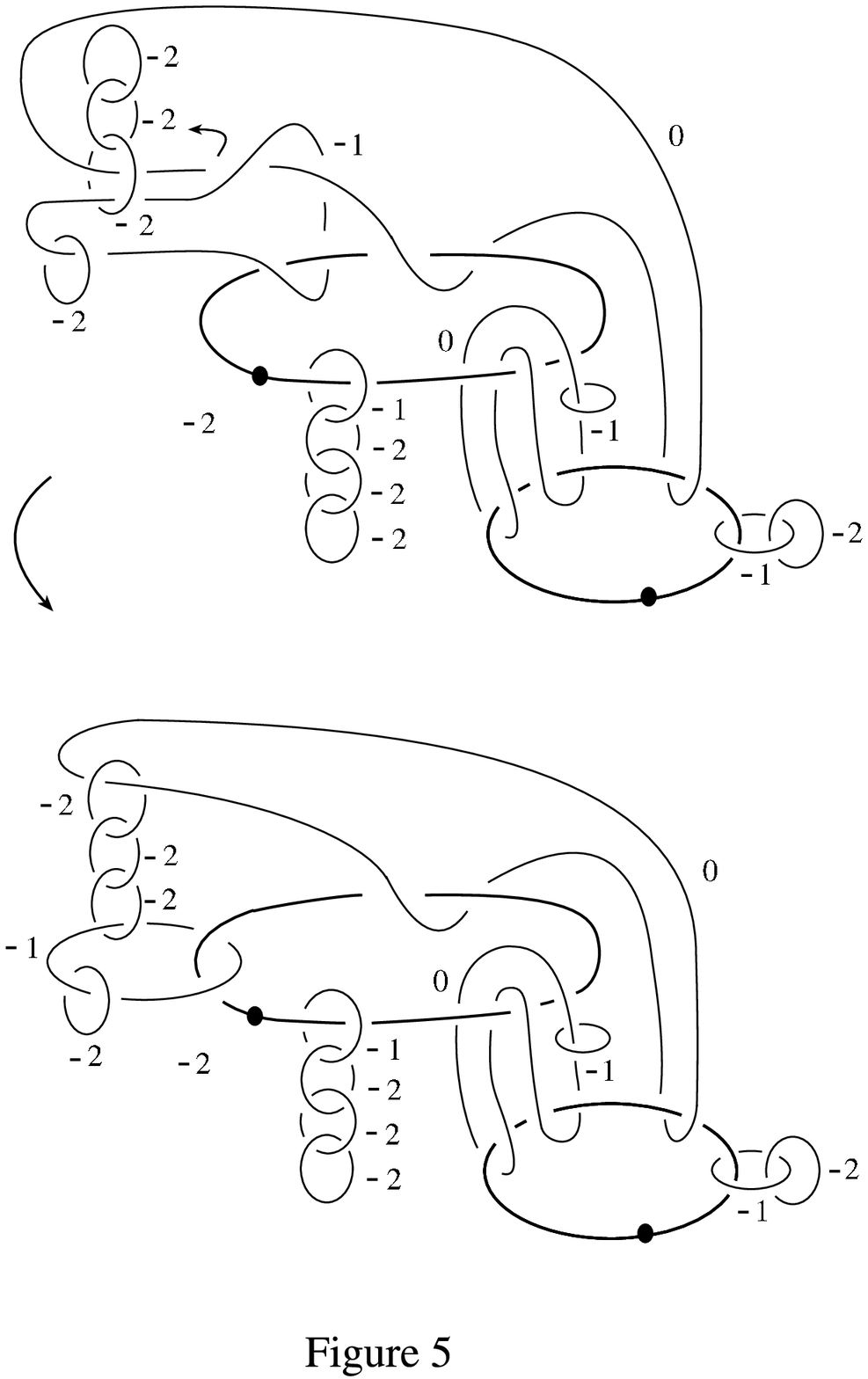}

\includegraphics[width=.95\textwidth]{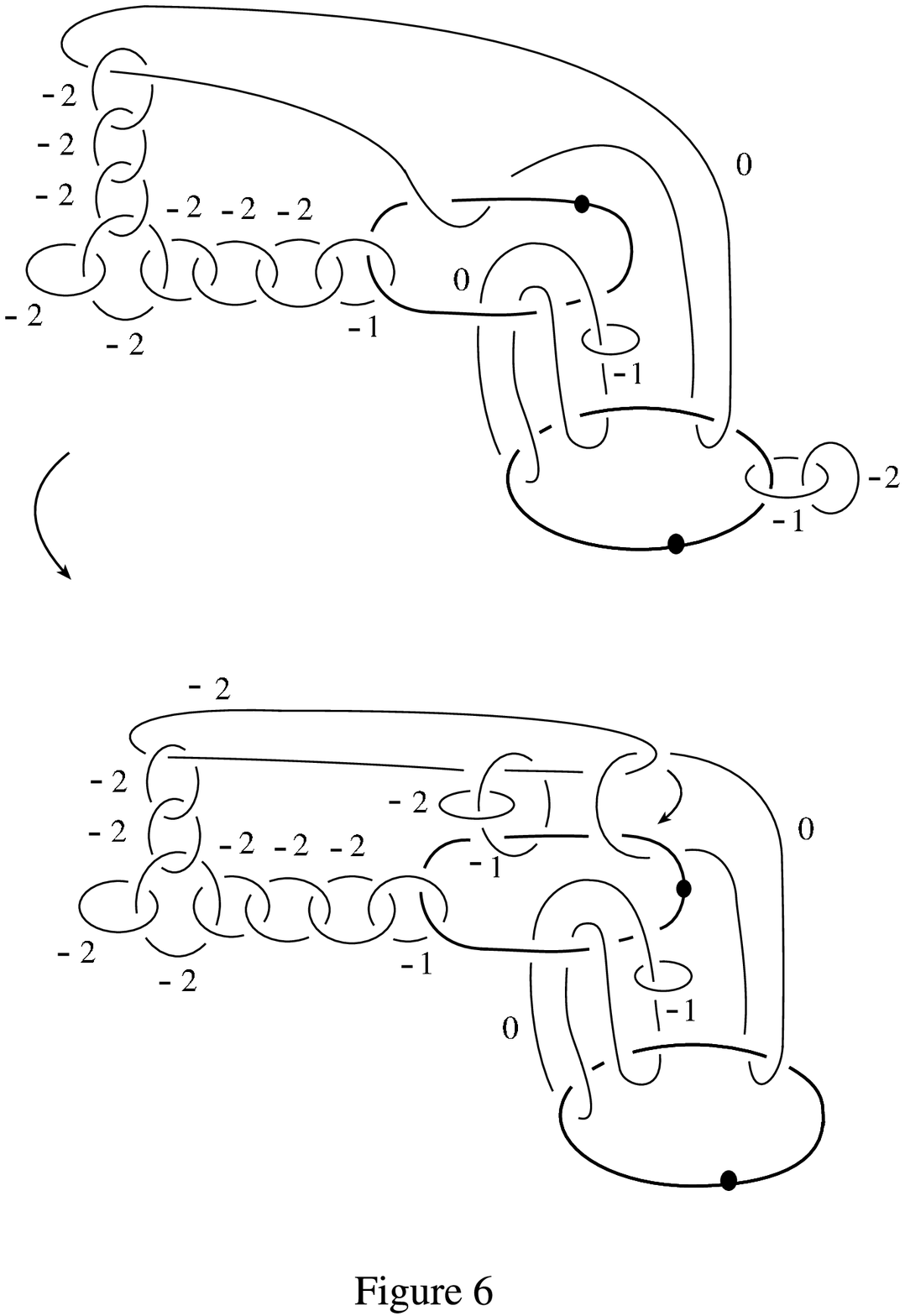}

\includegraphics{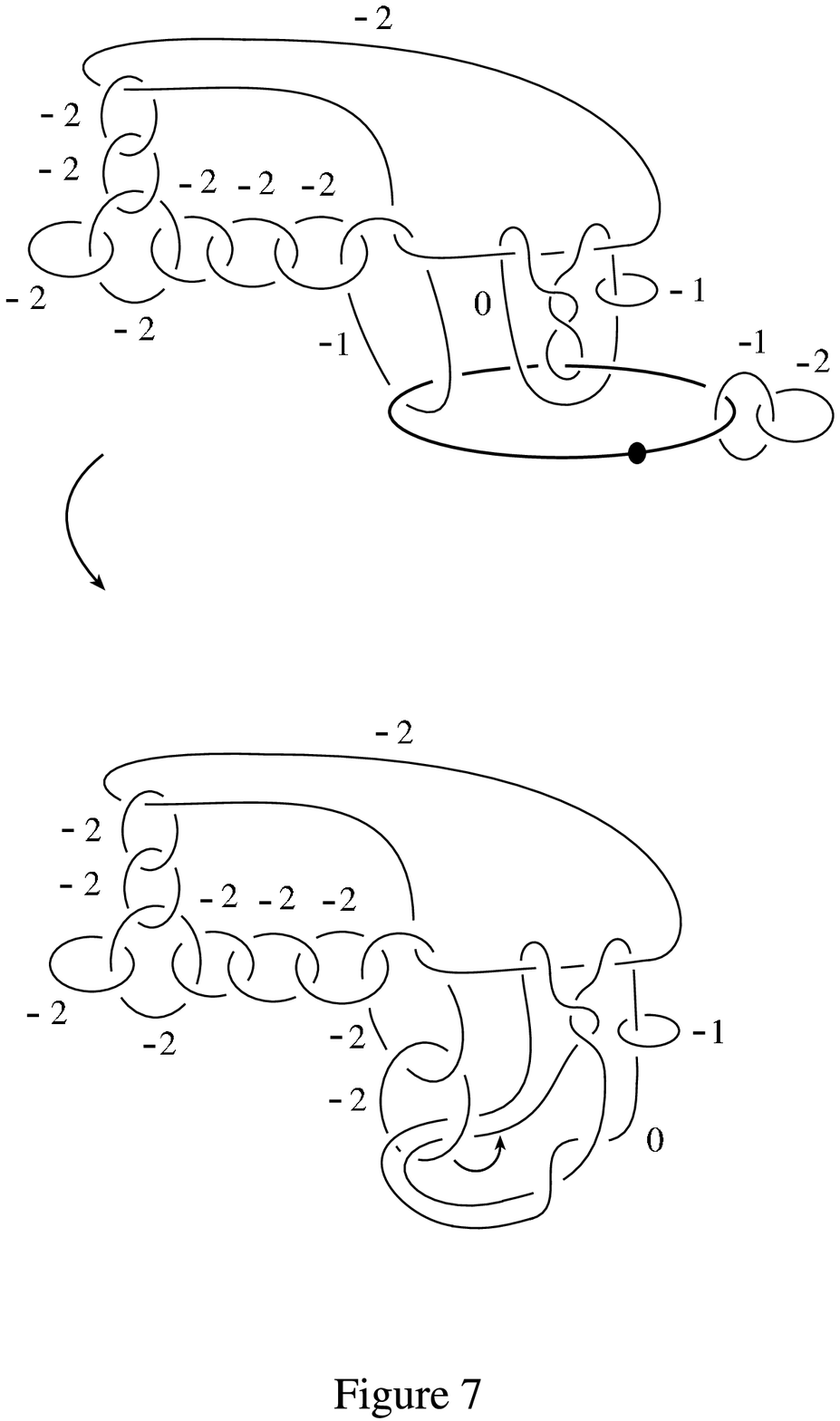}

\includegraphics{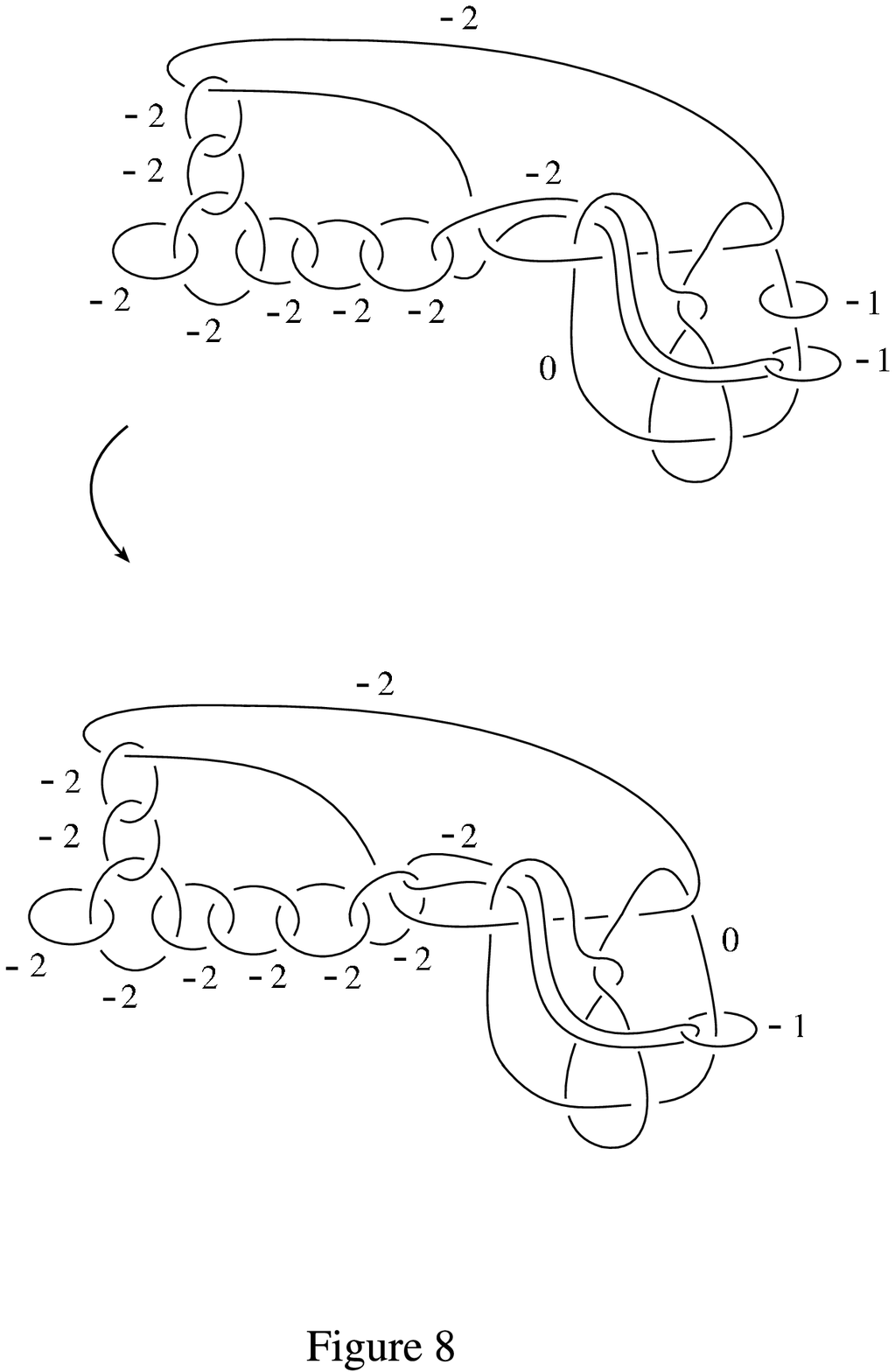}

\includegraphics{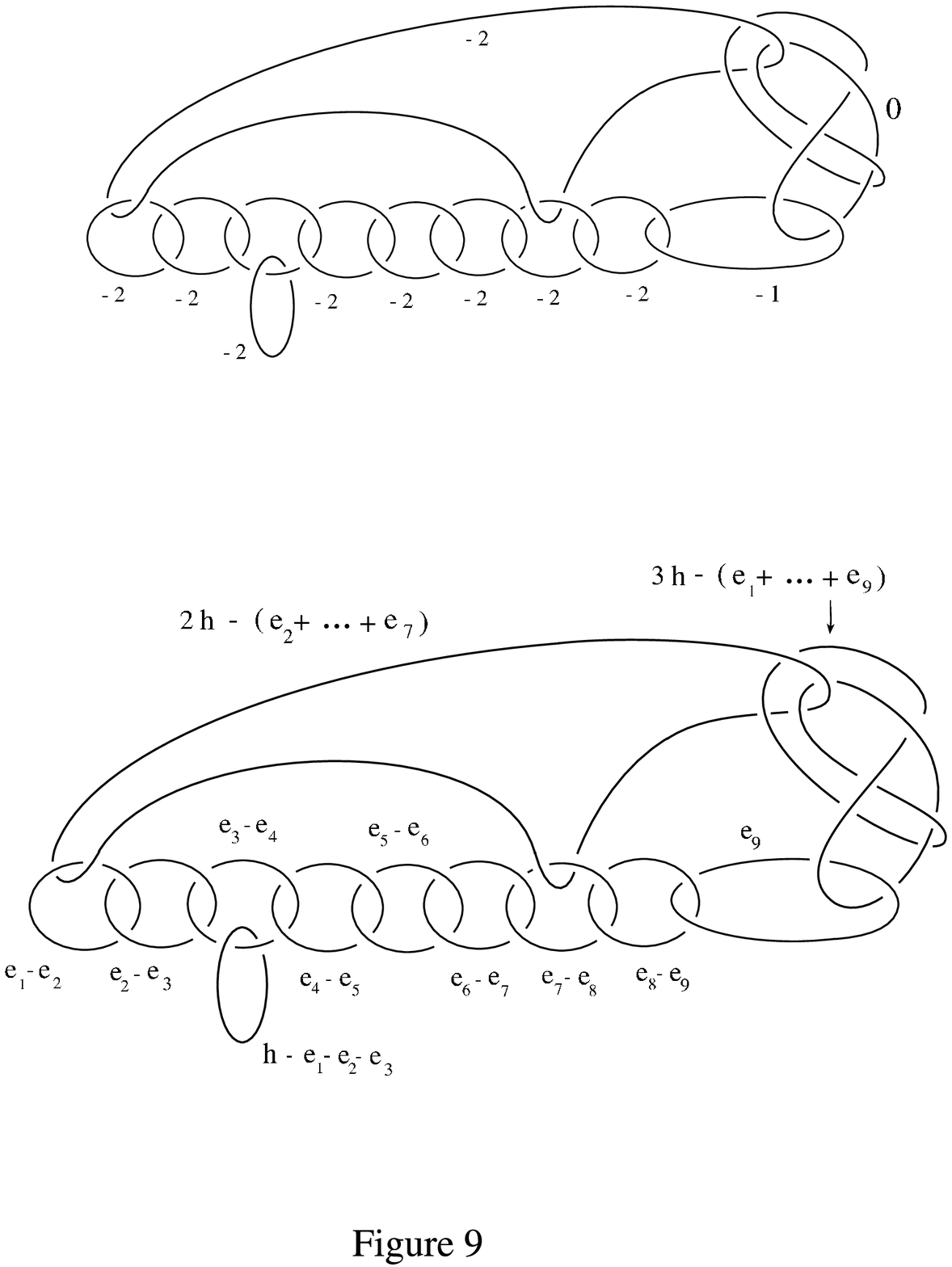}

\includegraphics[width=.95\textwidth]{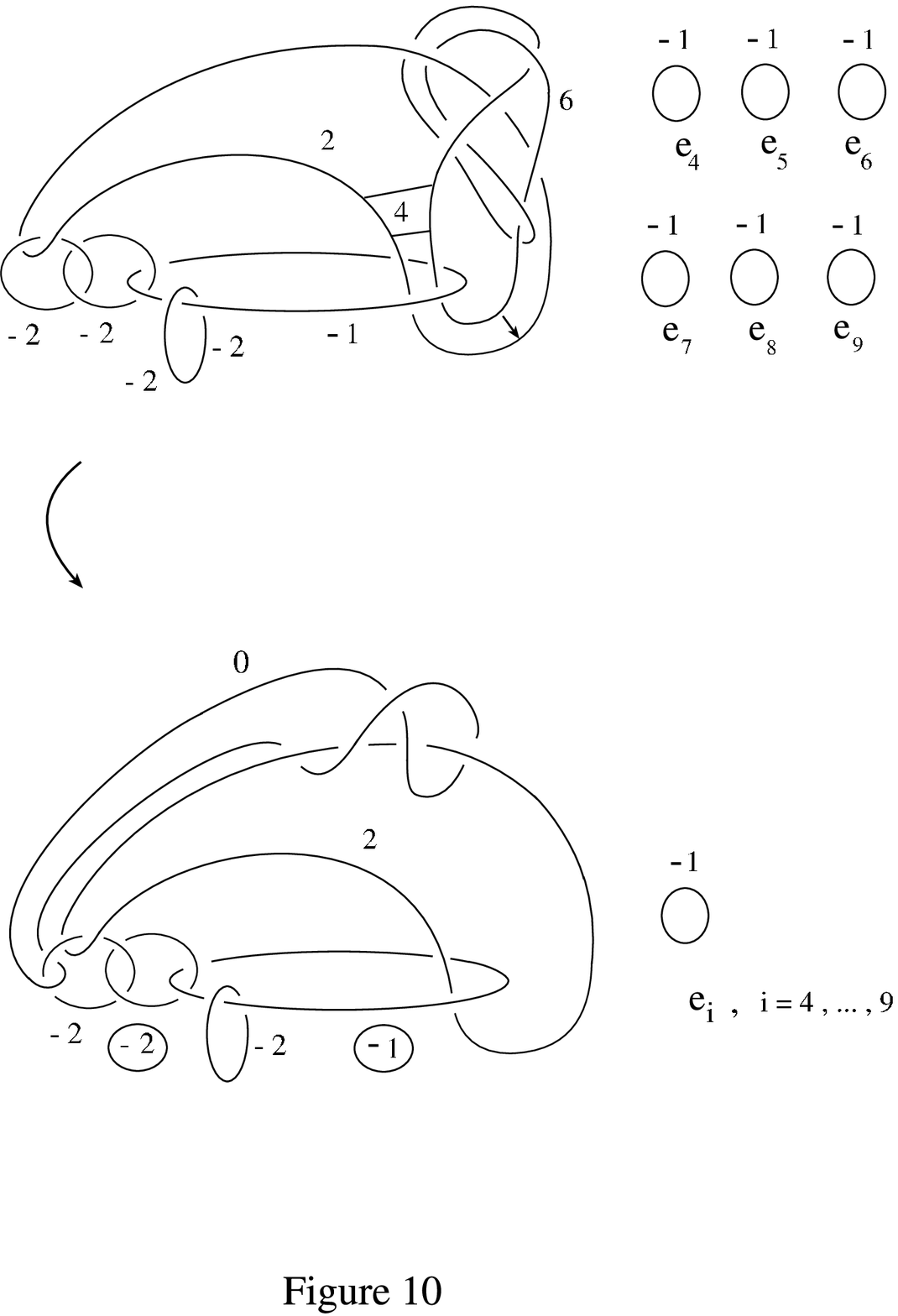}

\includegraphics[width=.95\textwidth]{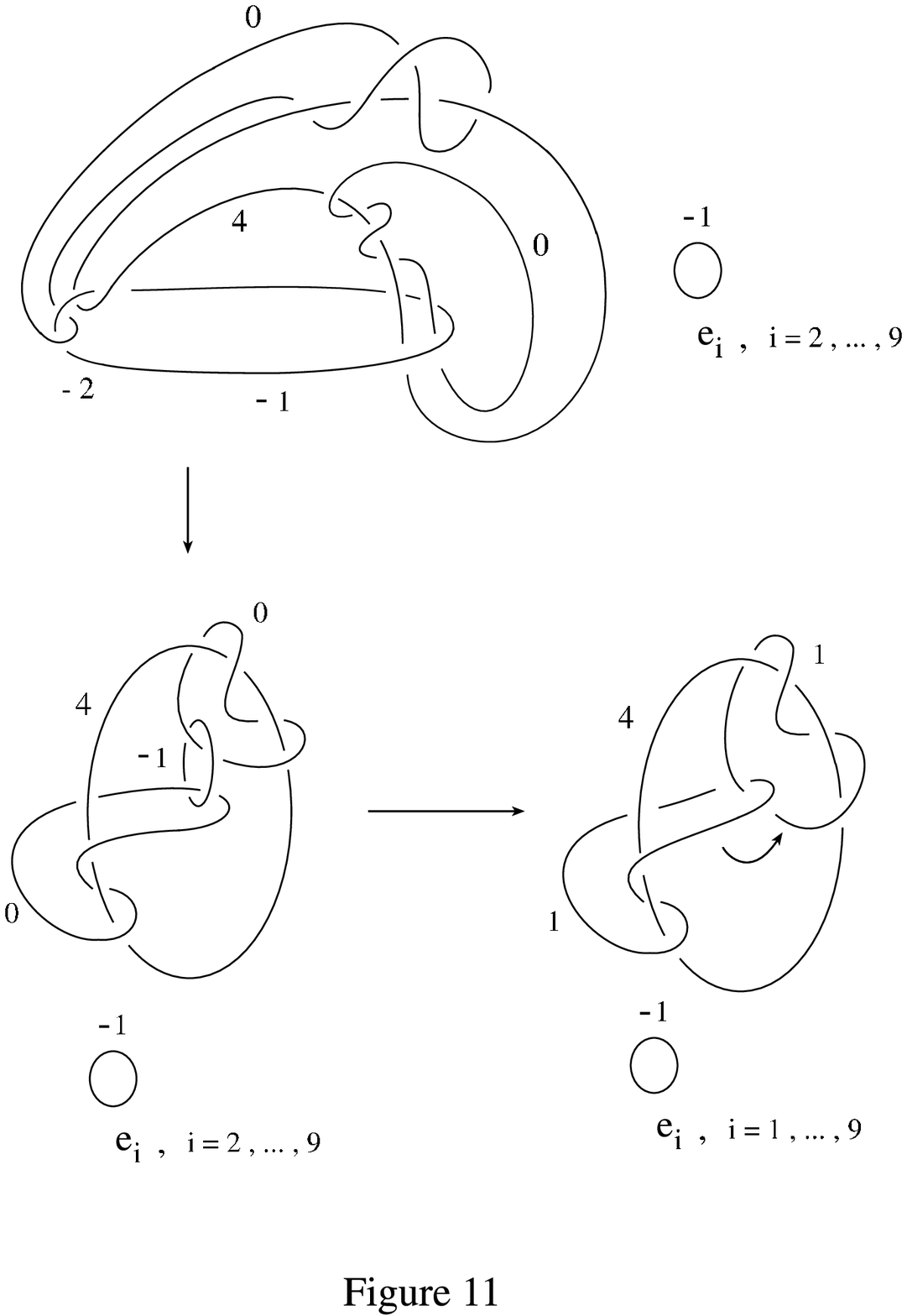}

\includegraphics[width=.9\textwidth]{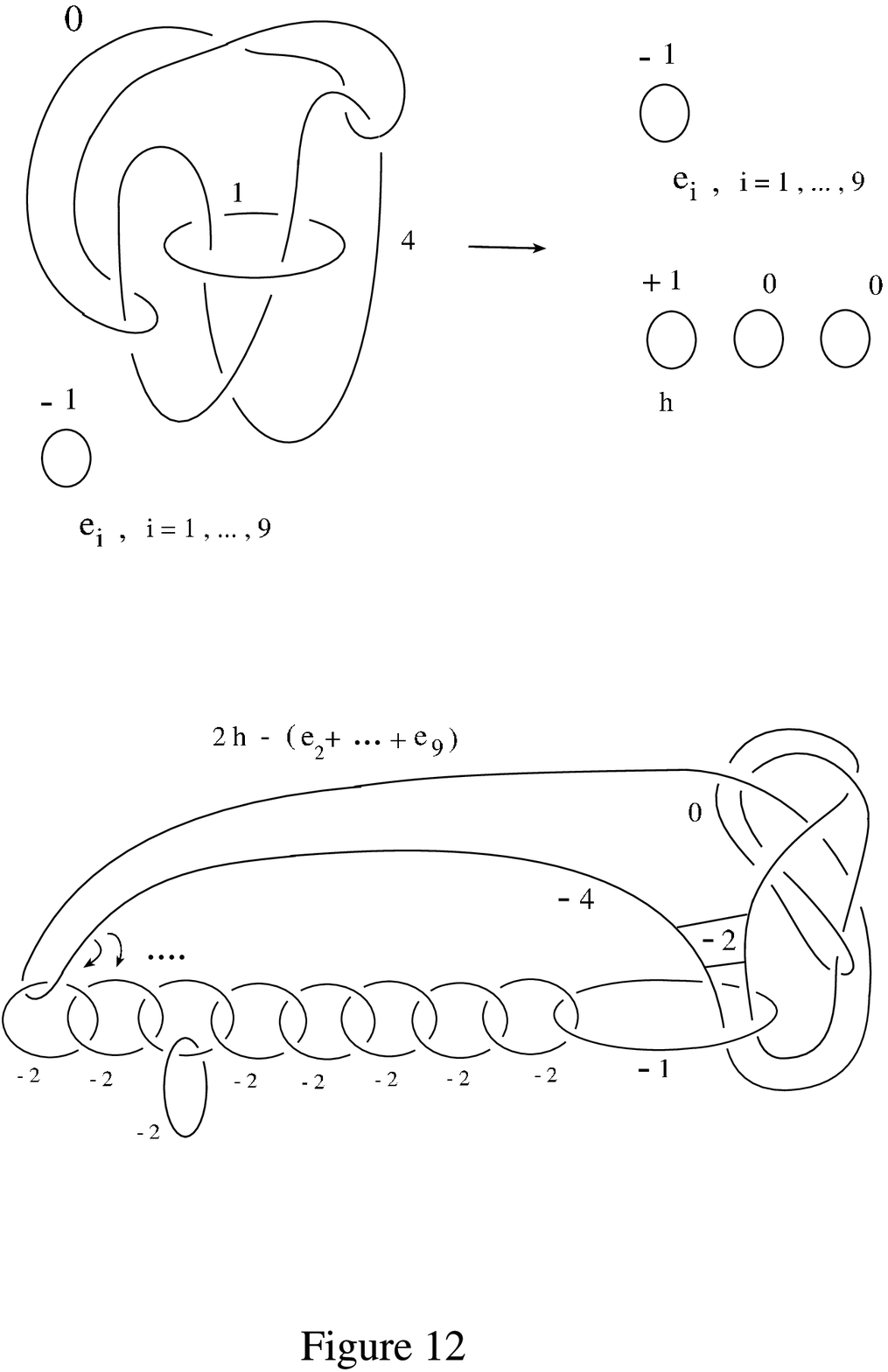}

\includegraphics{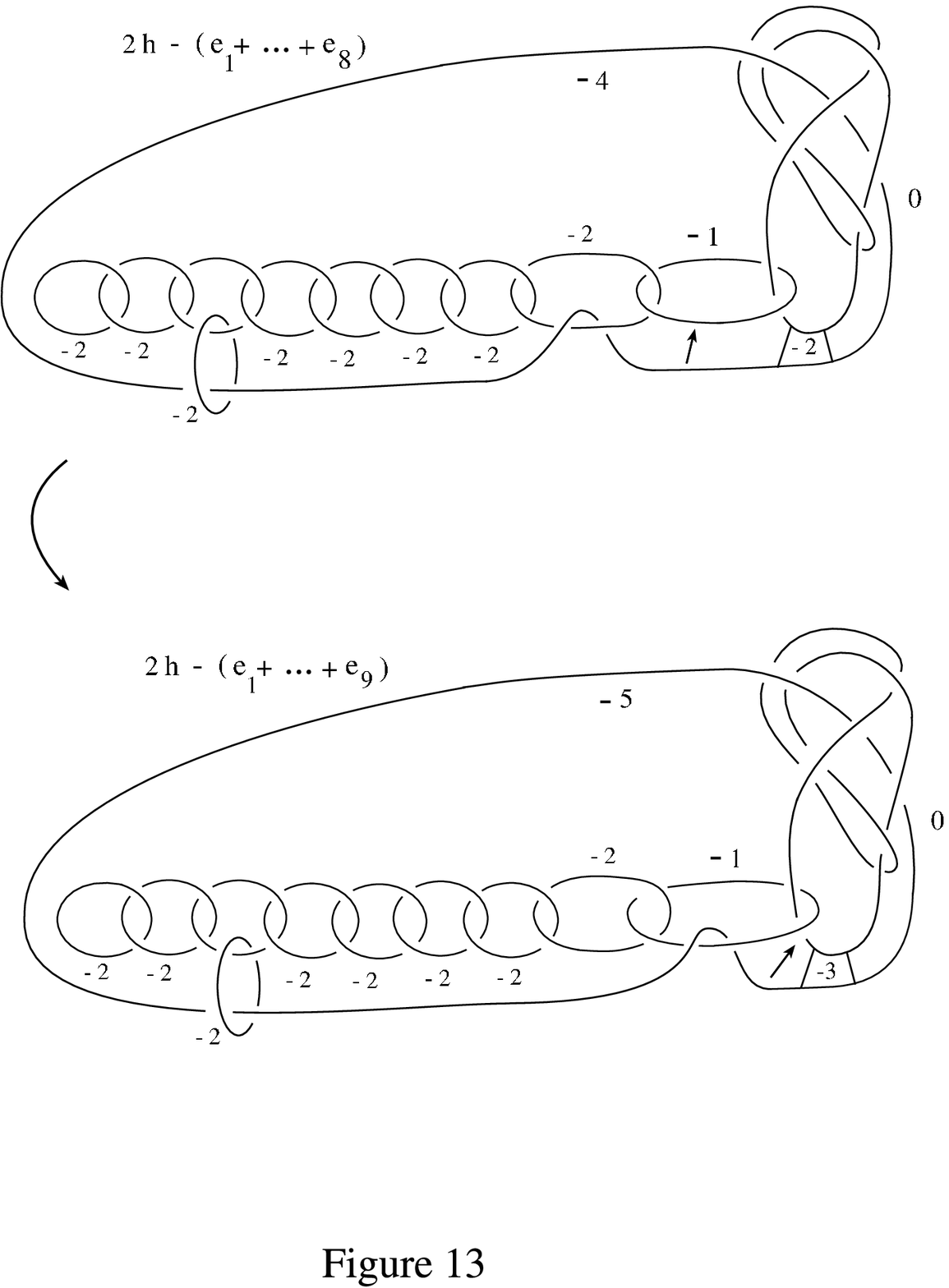}

\includegraphics{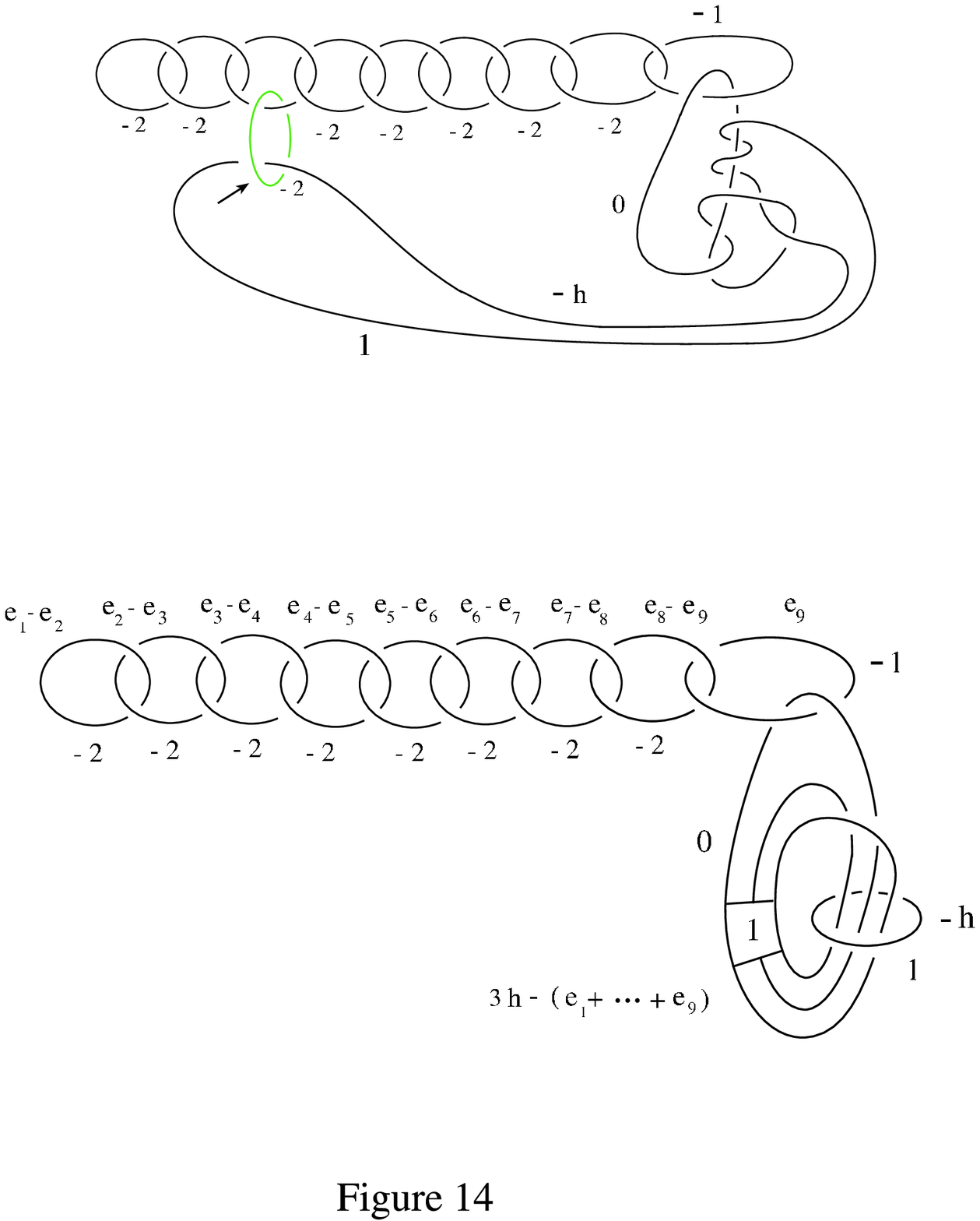}

\includegraphics[width=.85\textwidth]{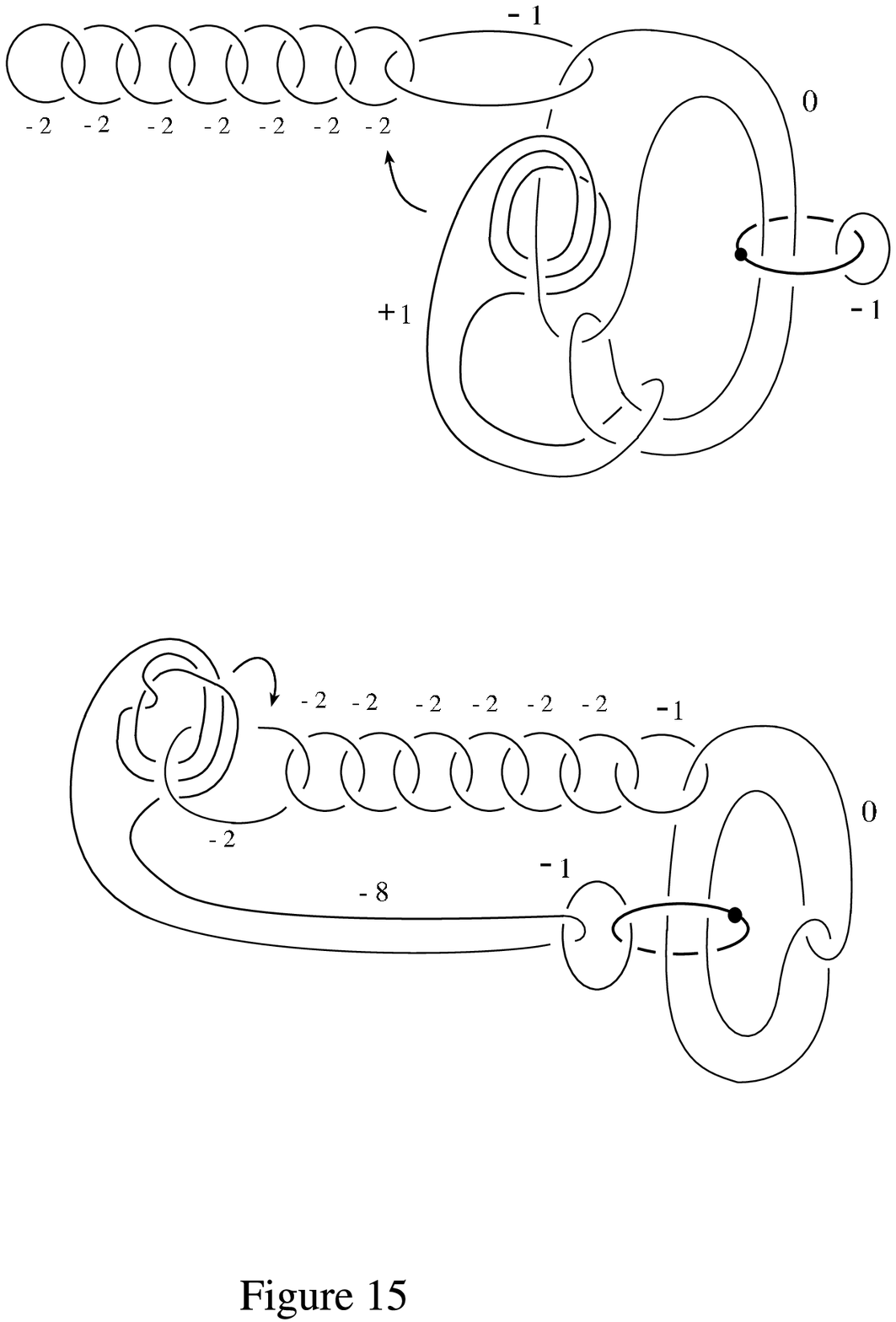}

\includegraphics[width=.8\textwidth]{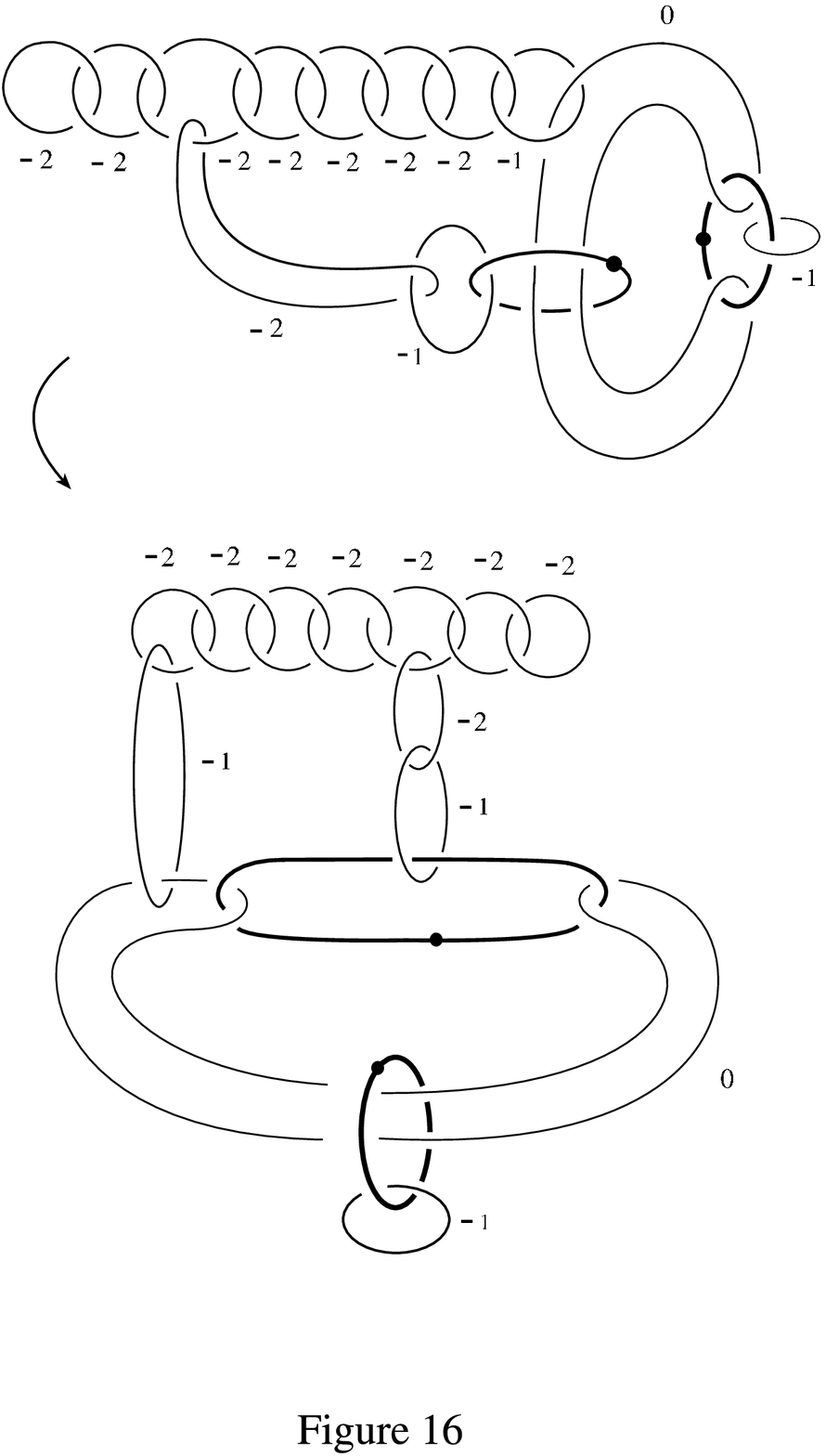}

\includegraphics[width=.75\textwidth]{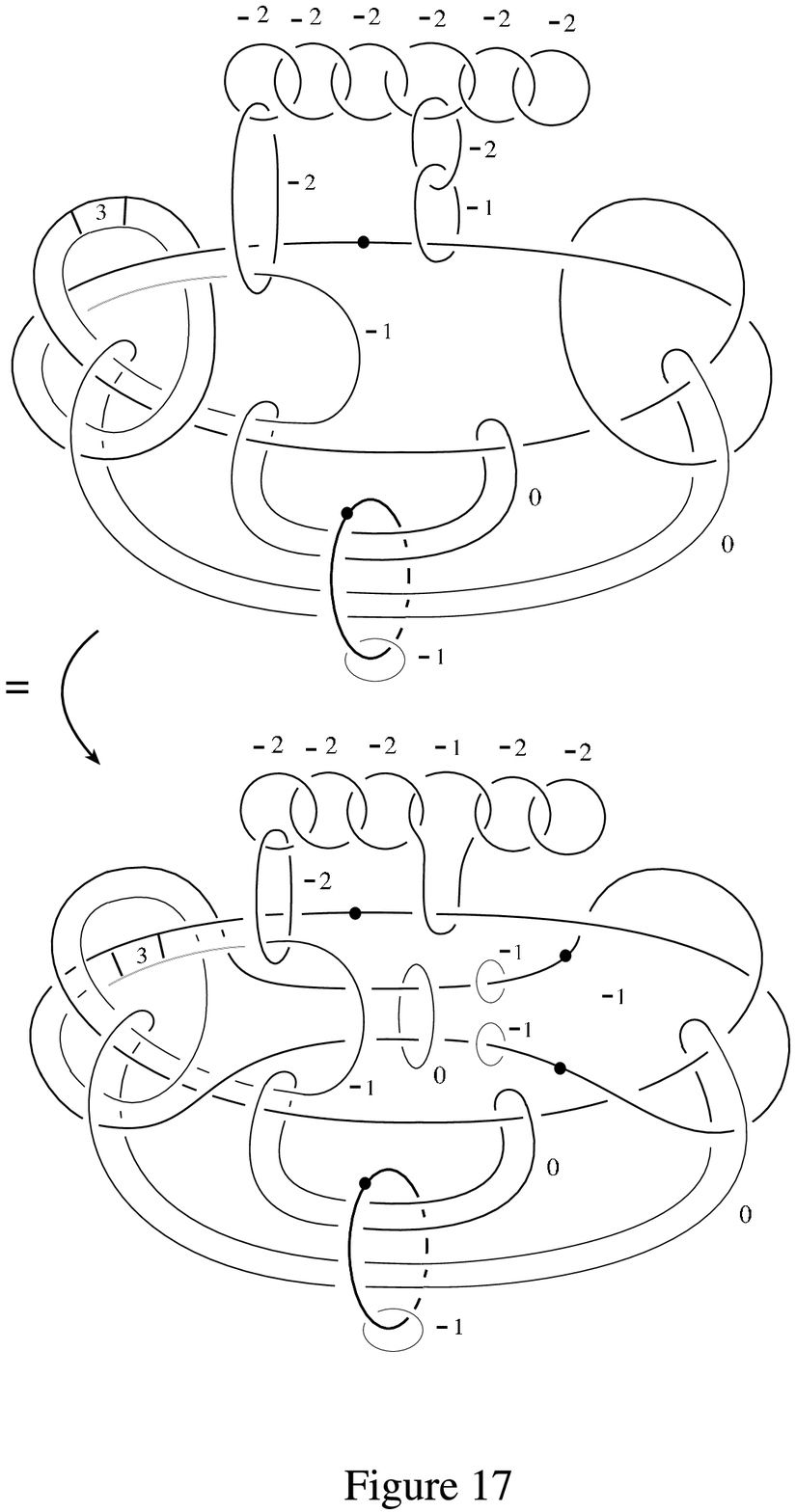}

\includegraphics[width=1\textwidth]{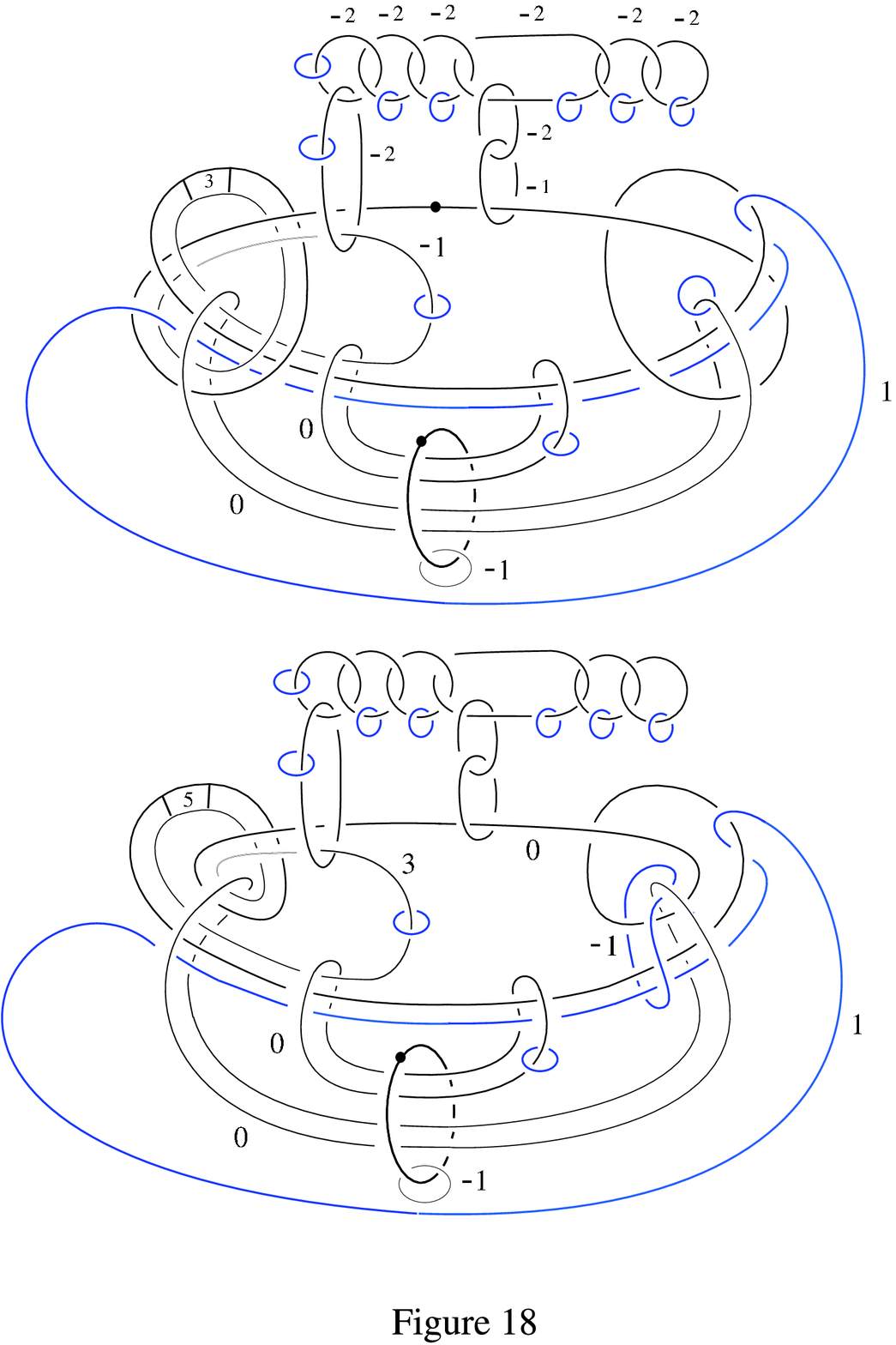}

\includegraphics[width=1\textwidth]{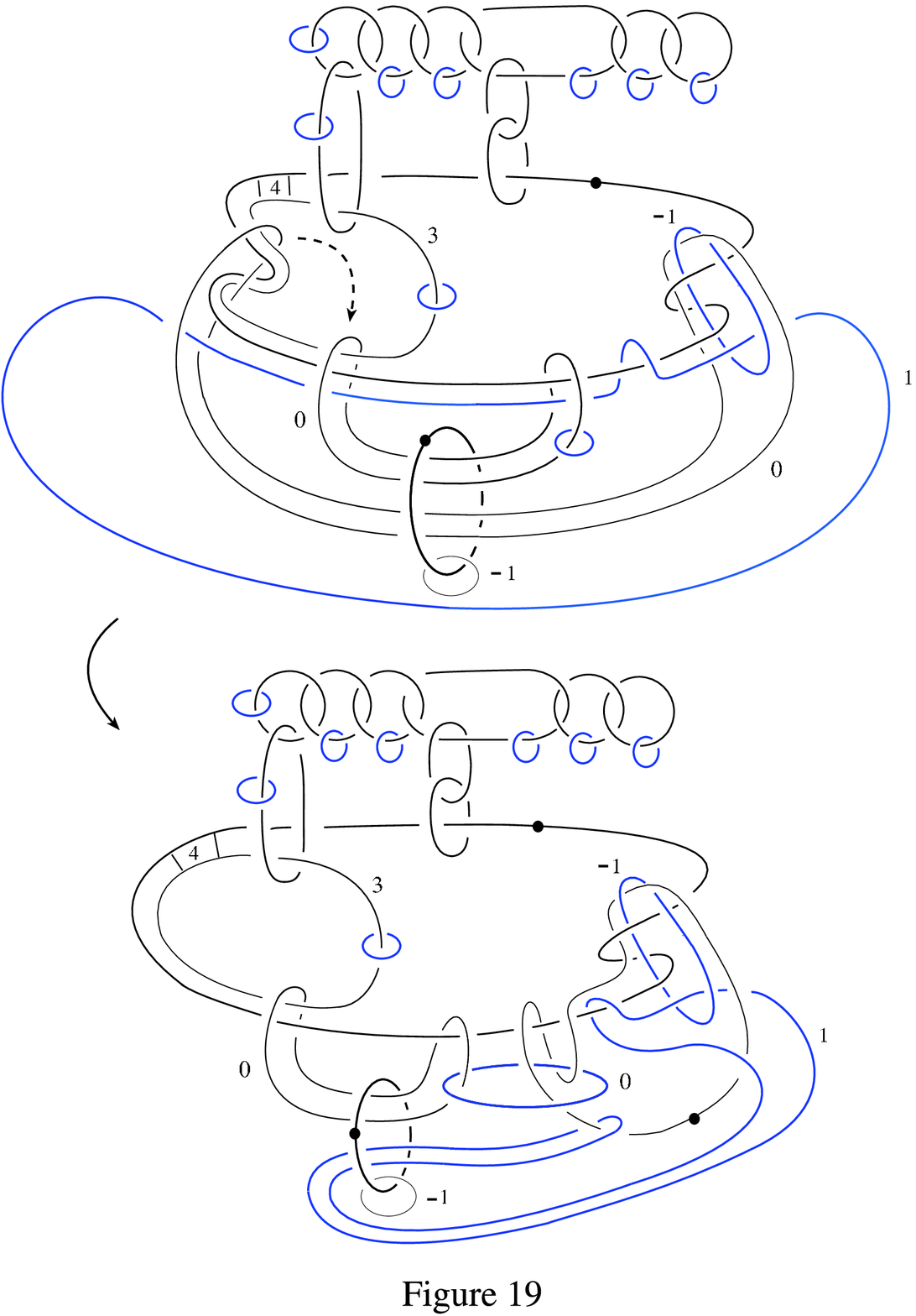}

\includegraphics{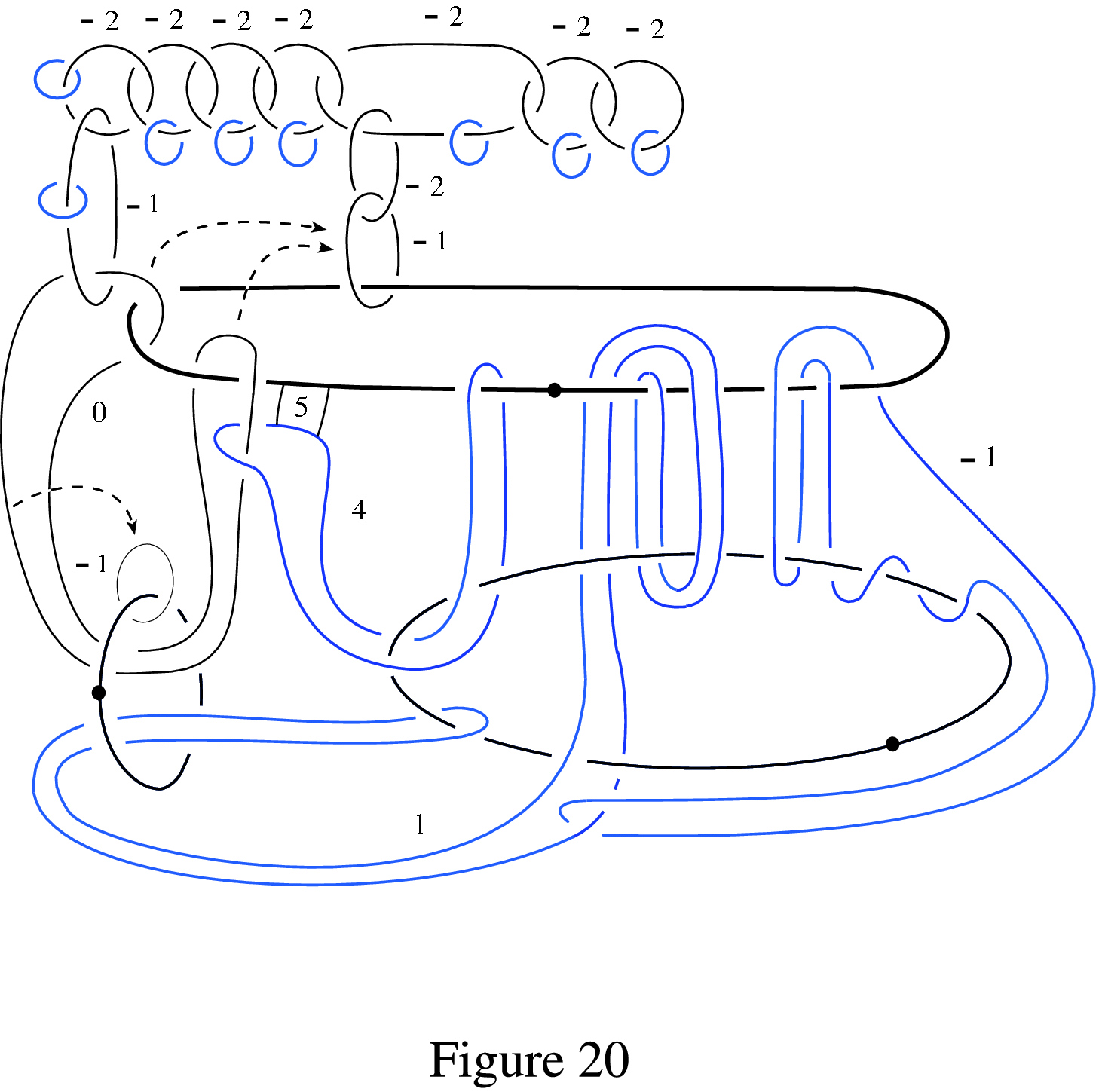}

\includegraphics{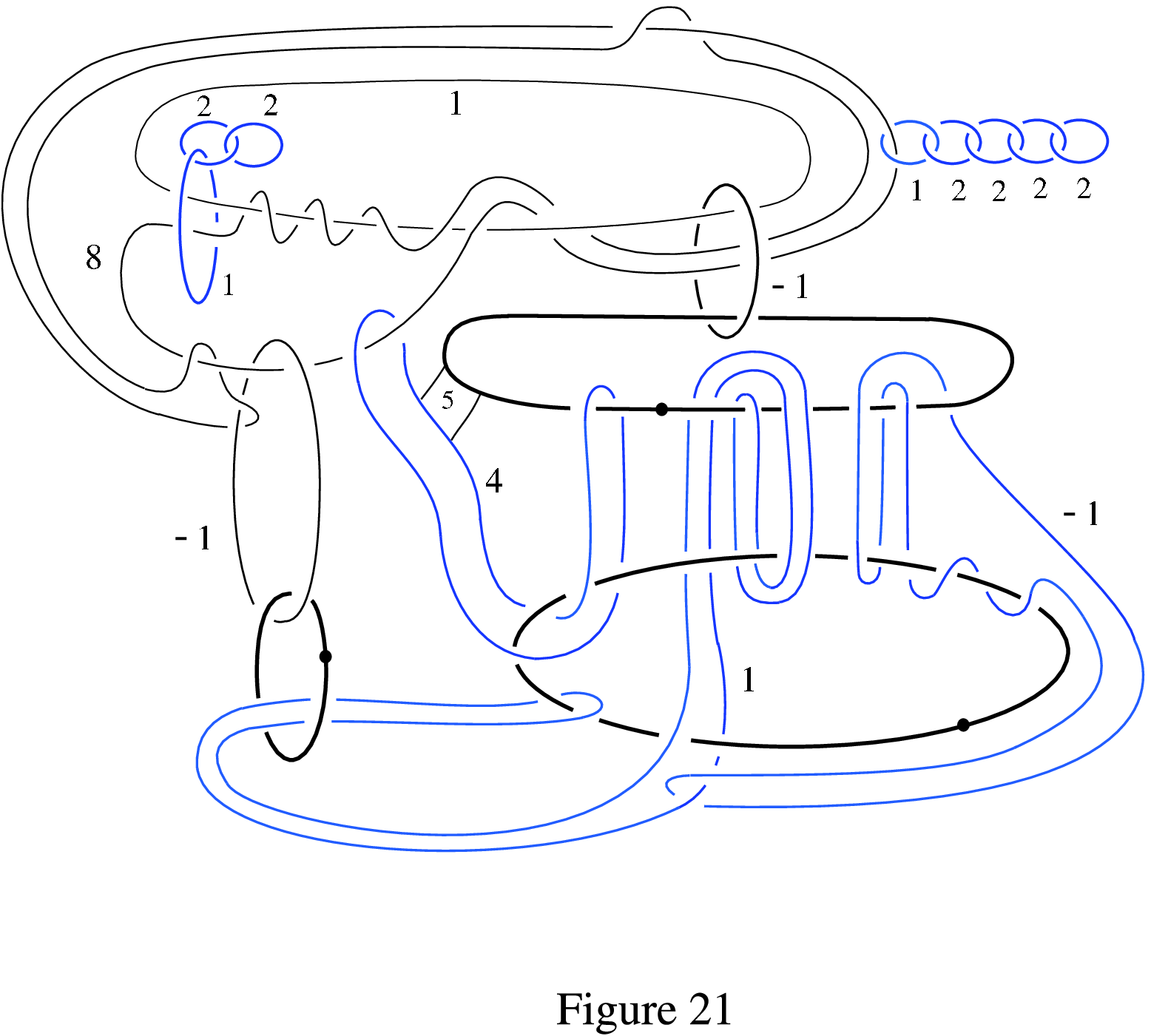}

\includegraphics{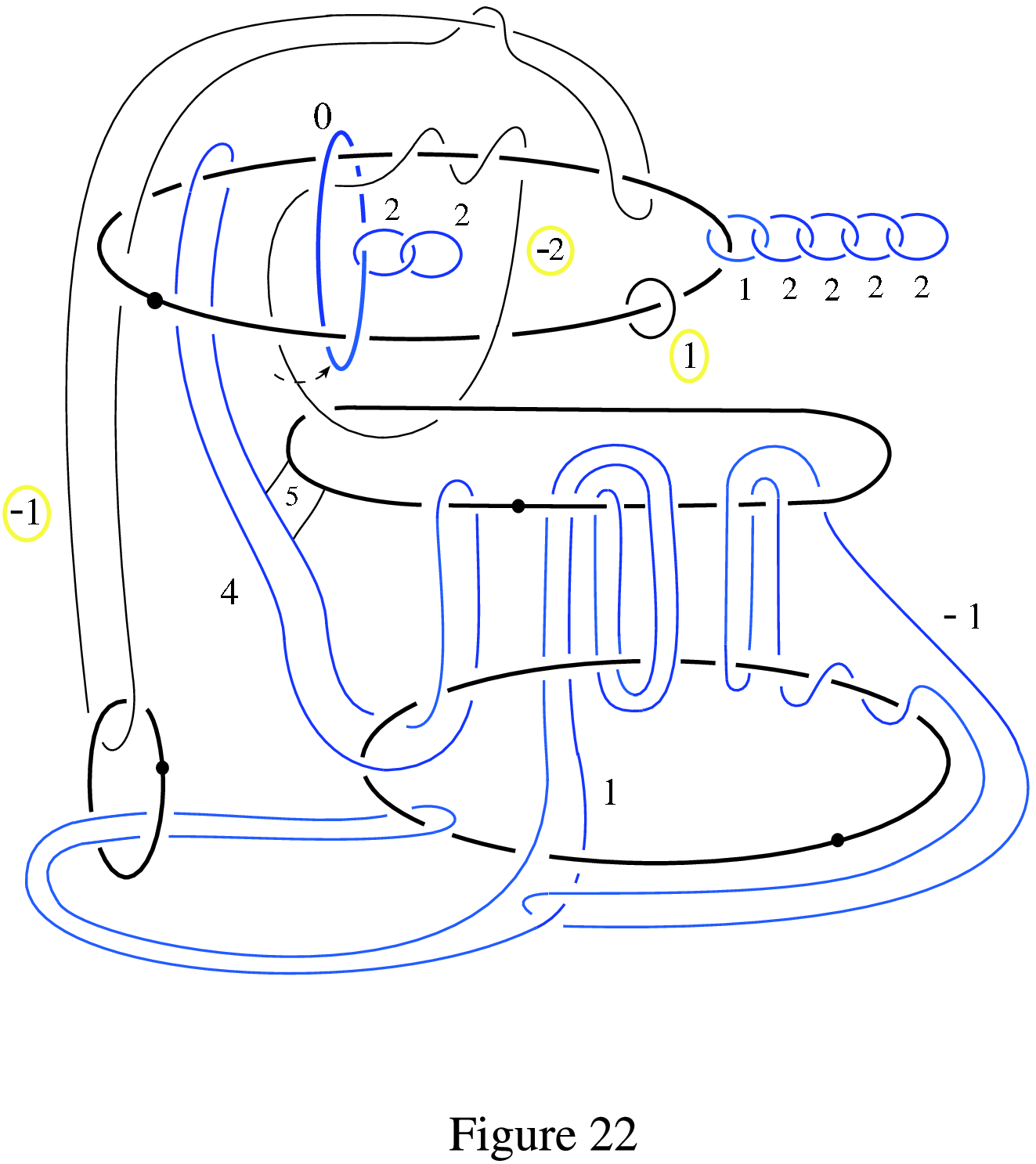}

\includegraphics{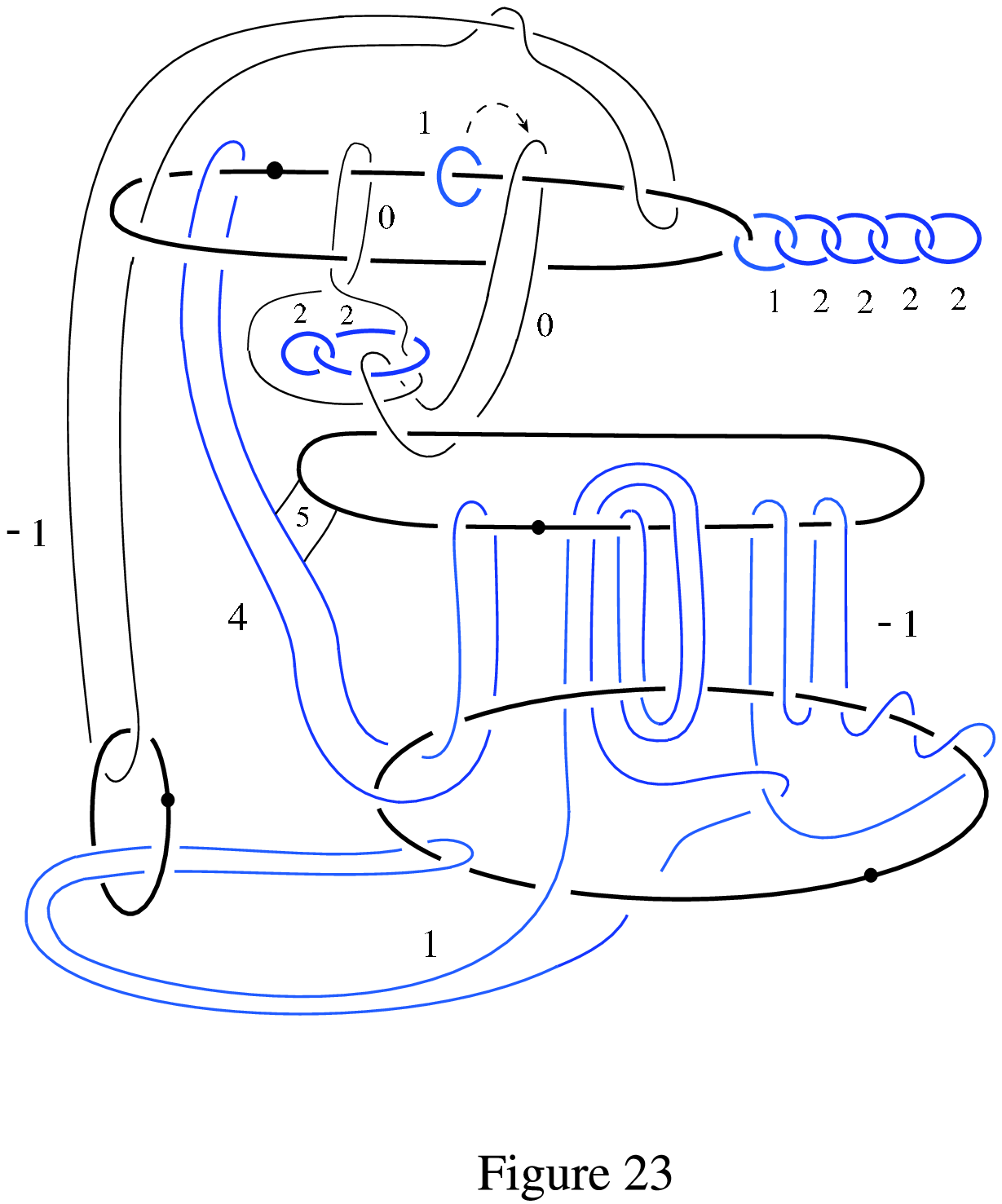}

\includegraphics{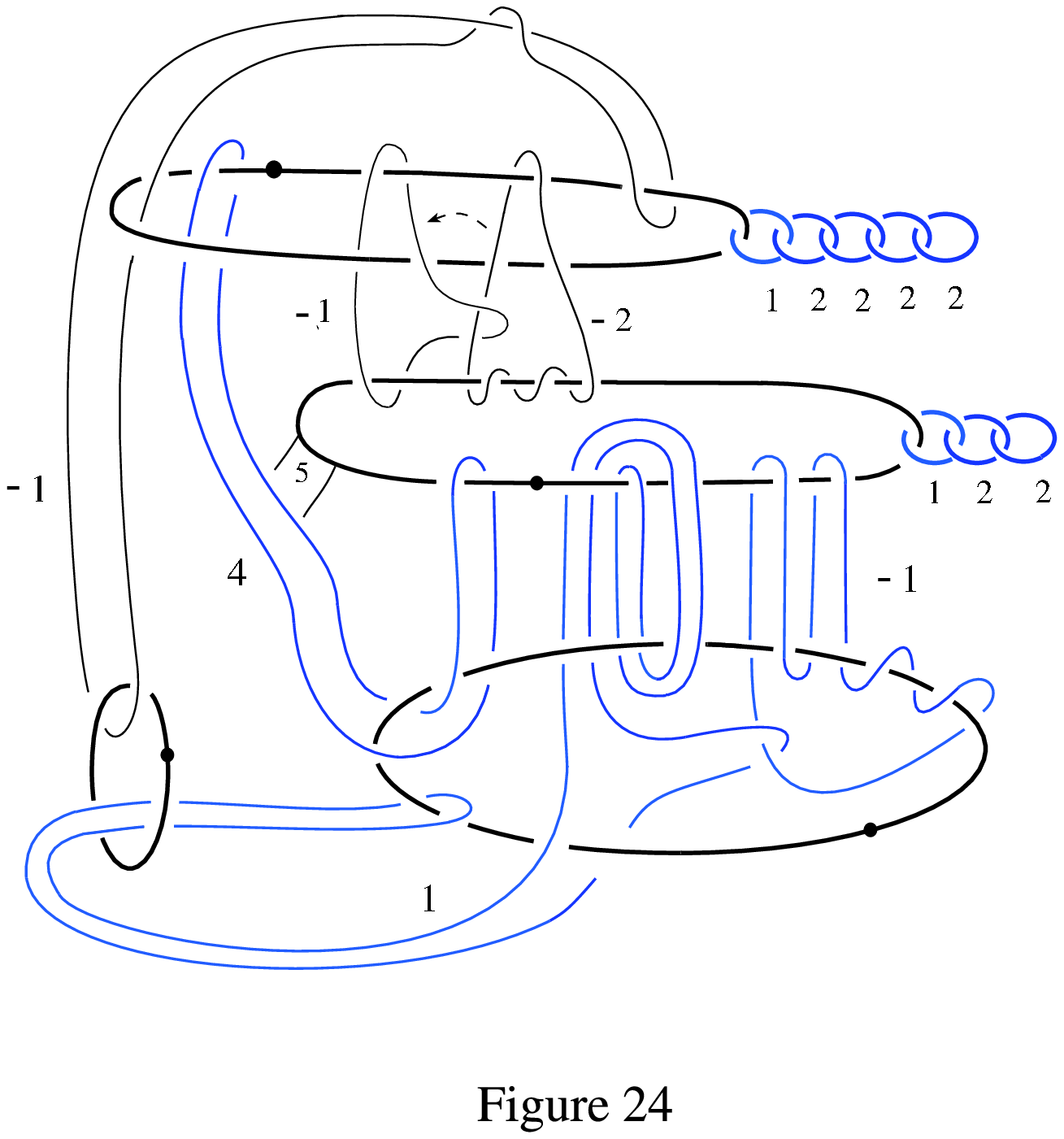}

\includegraphics{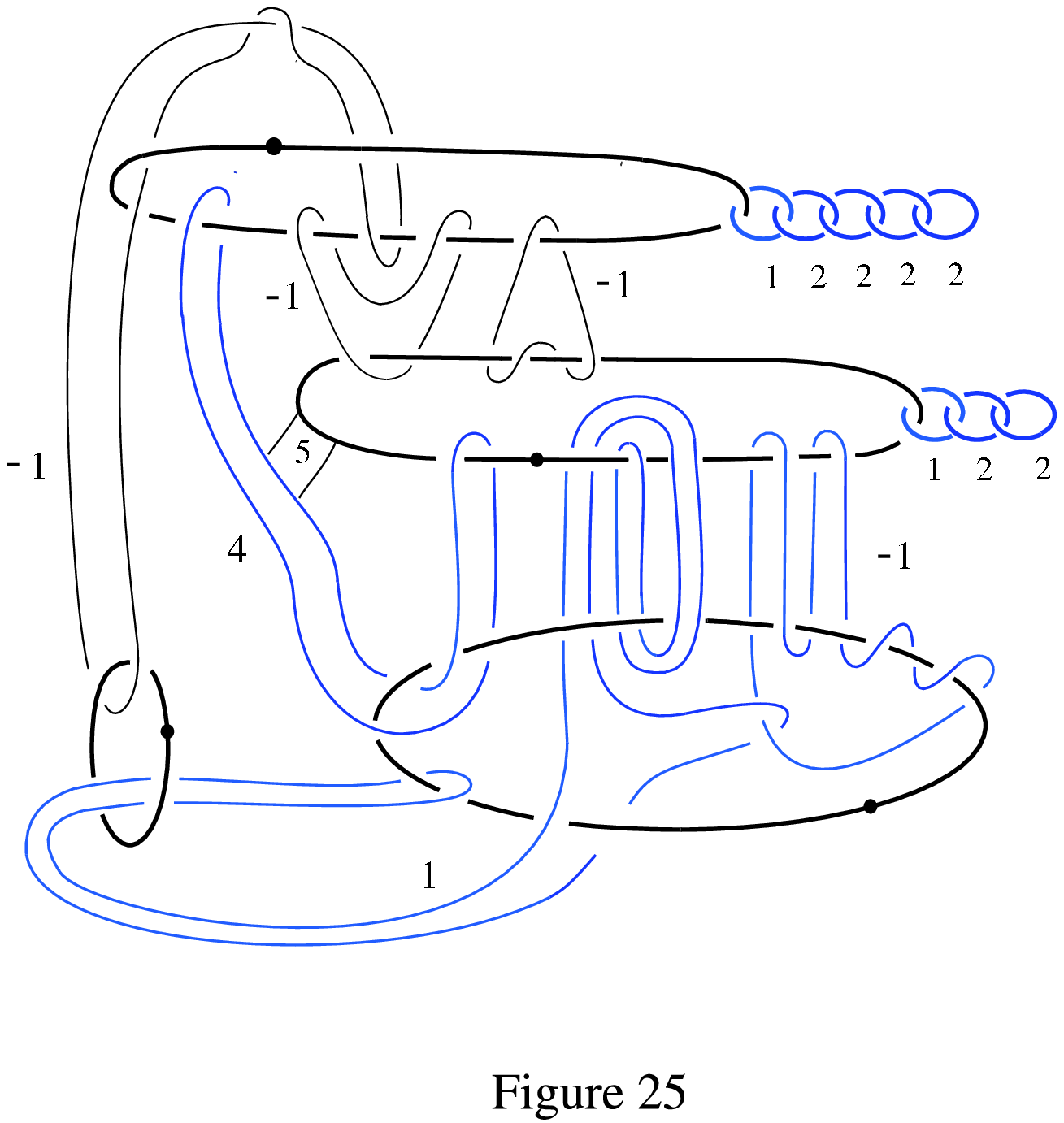}

\includegraphics[width=.95\textwidth]{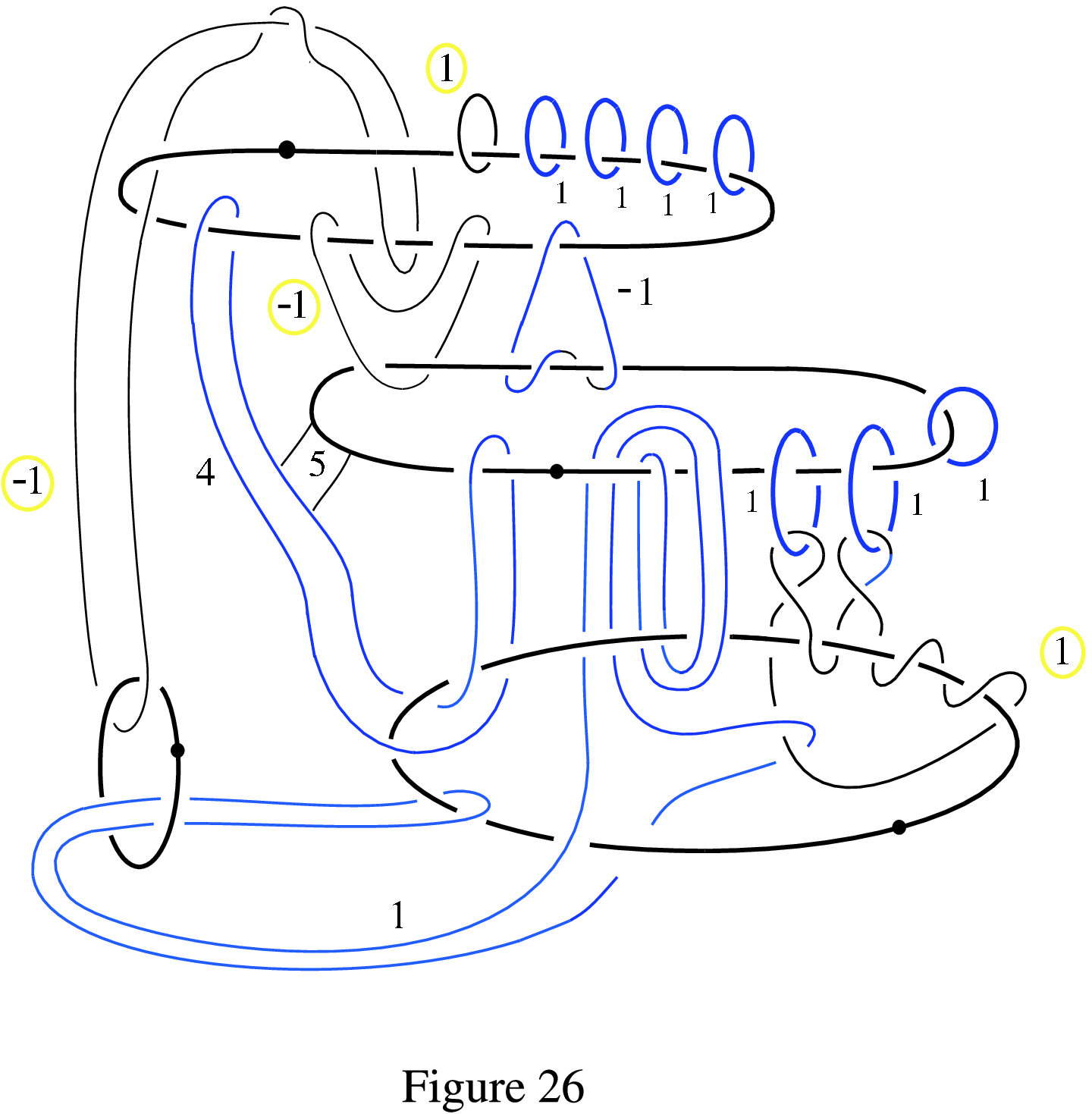}

\includegraphics{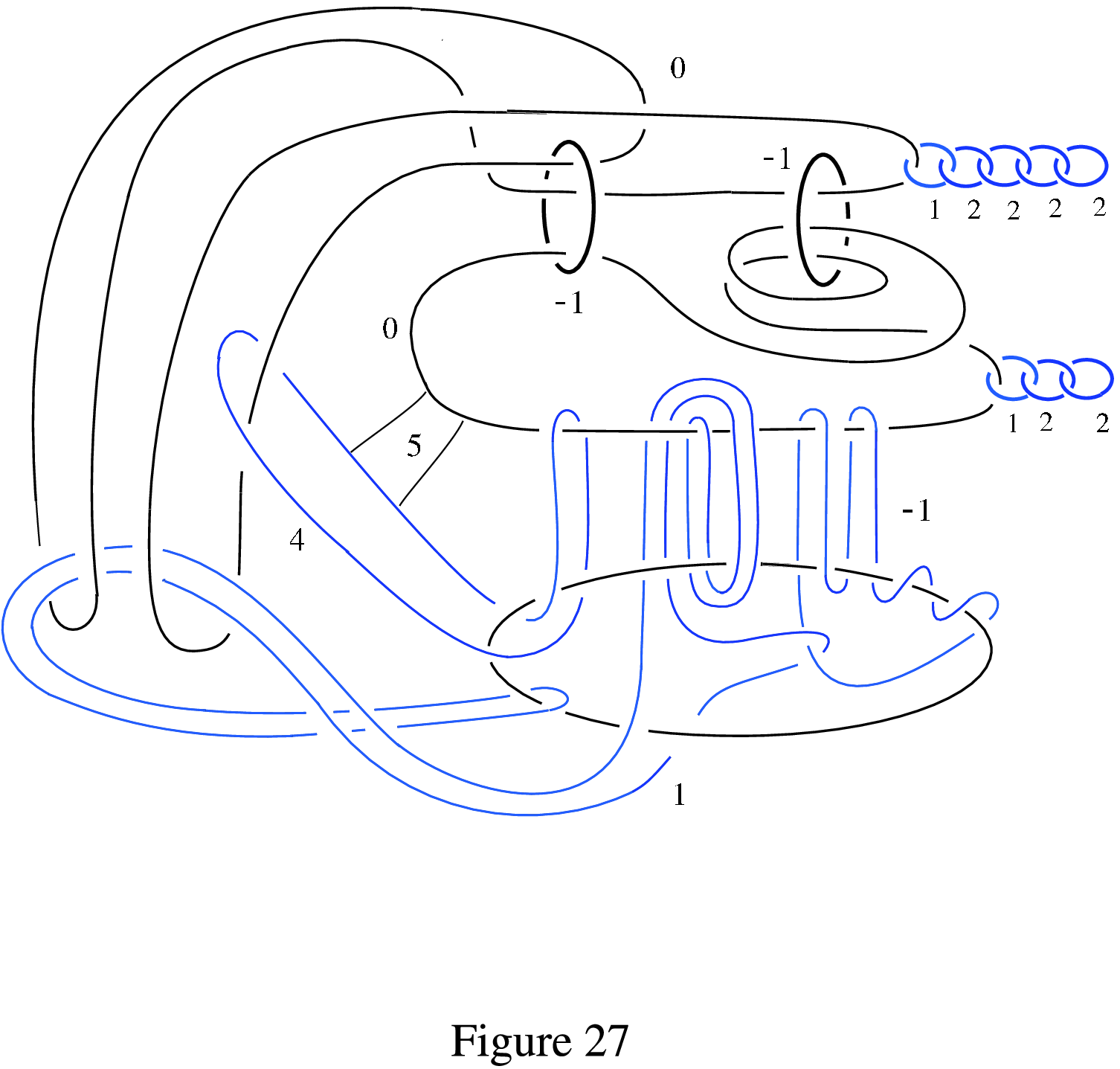}

\includegraphics{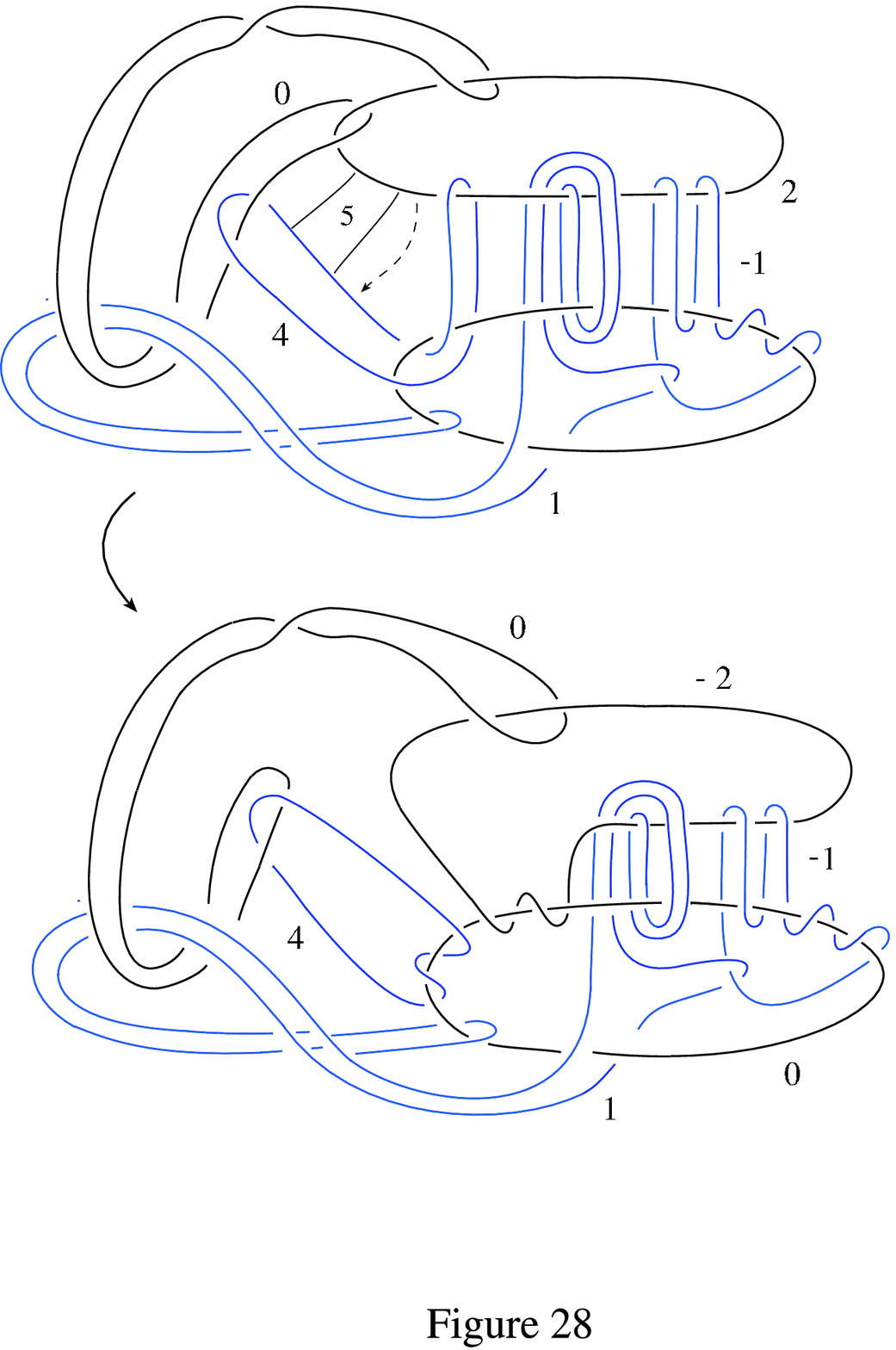}

\includegraphics[width=.8\textwidth]{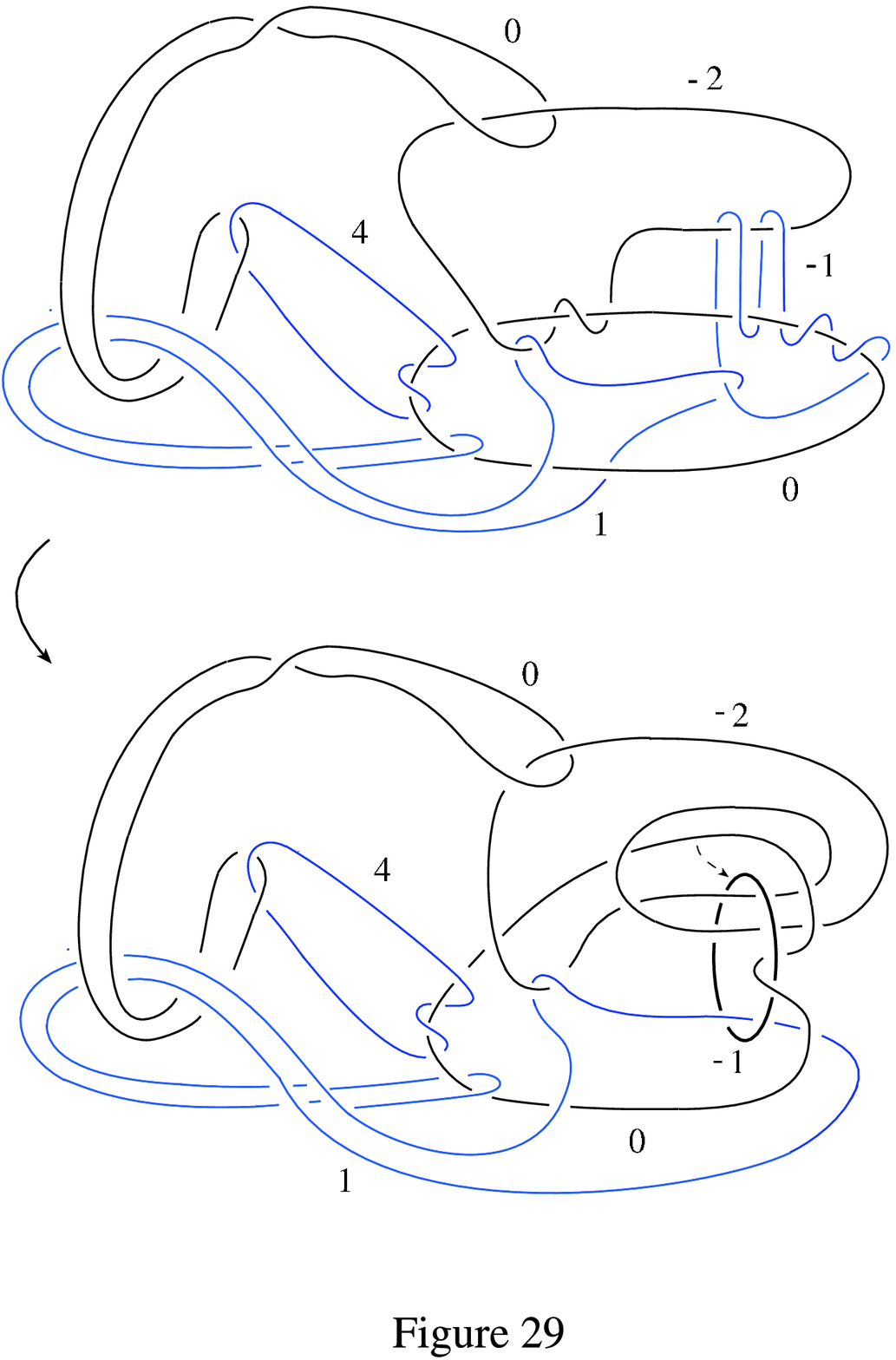}

\includegraphics[width=.8\textwidth]{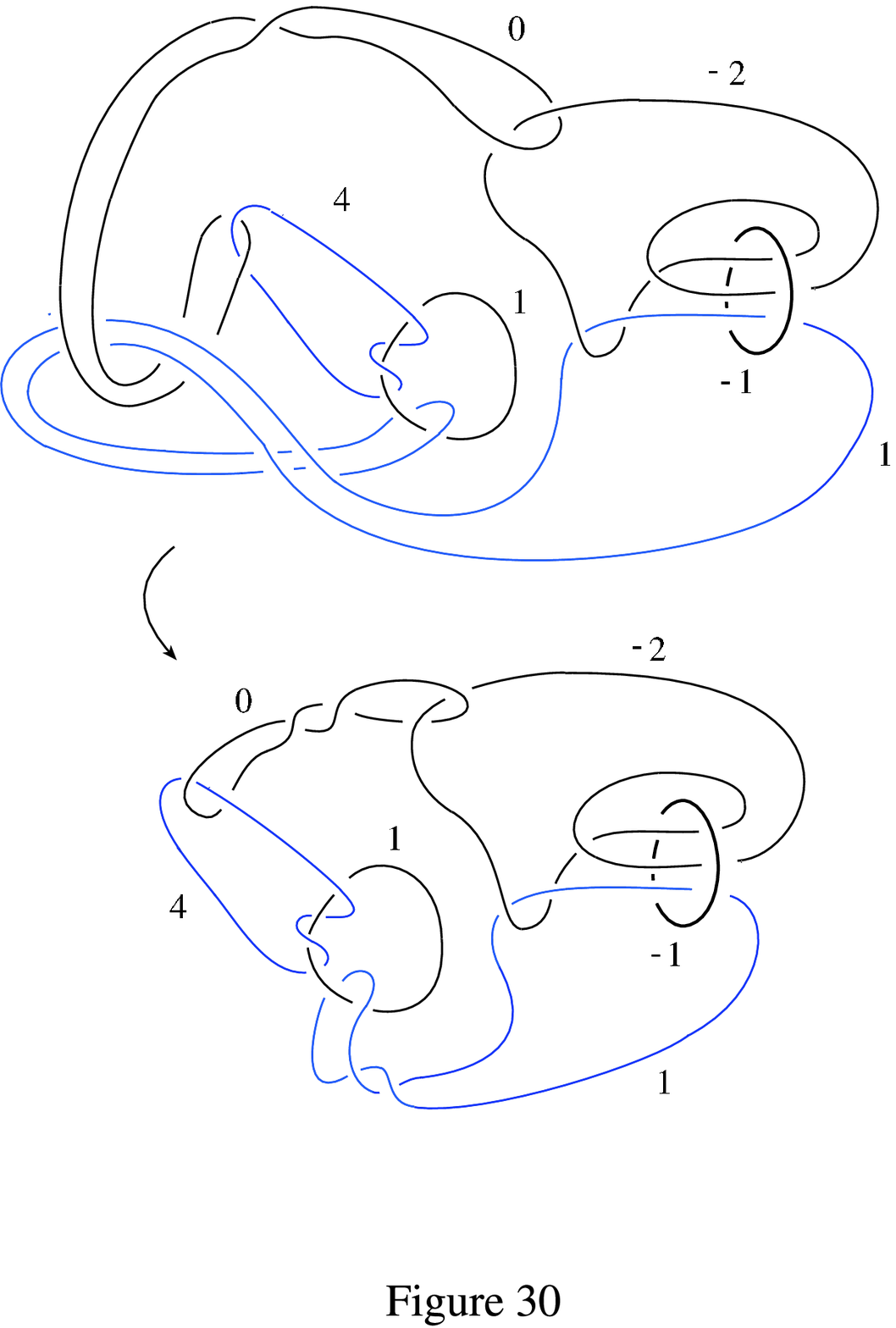}

\includegraphics[width=.9\textwidth]{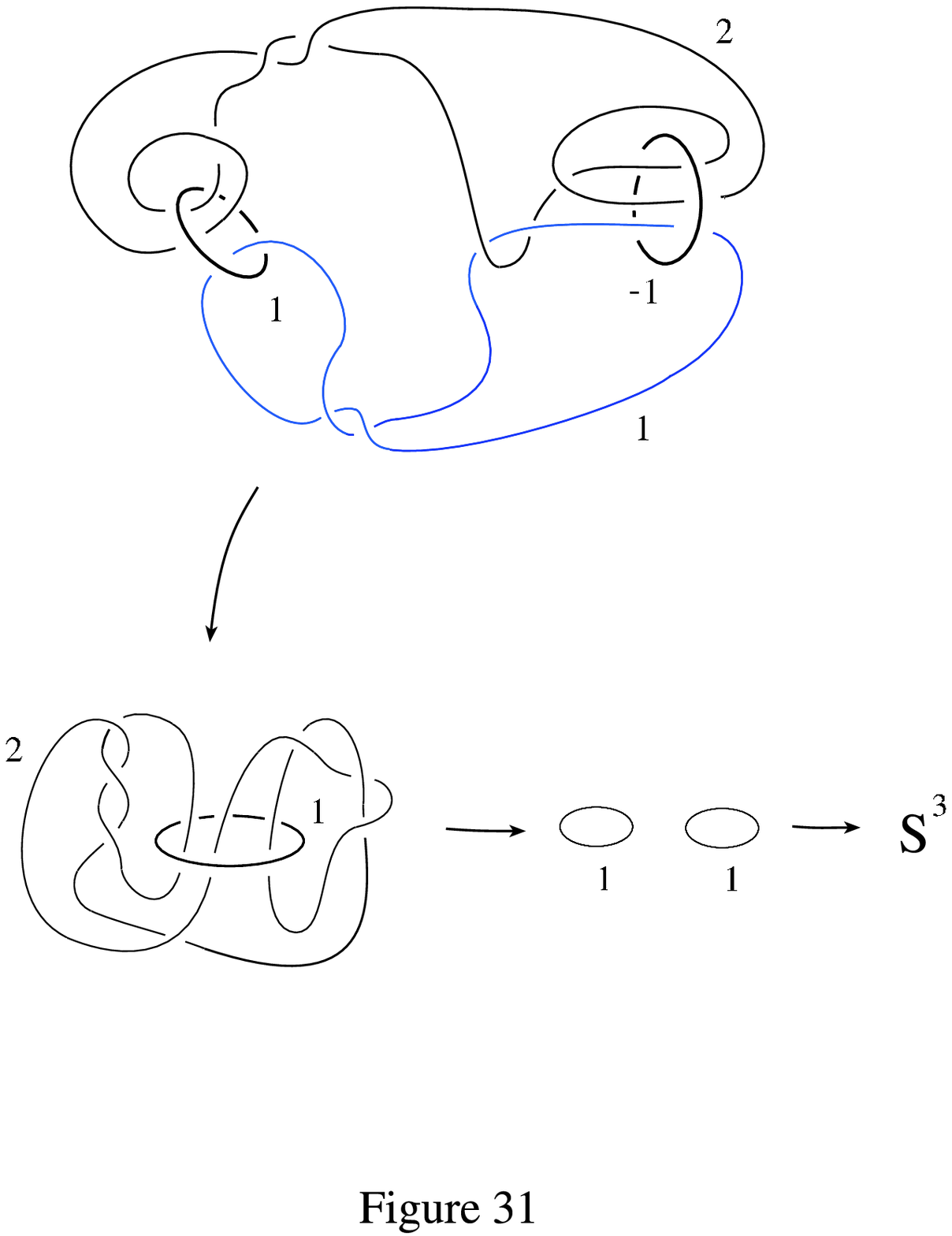}

\includegraphics[width=.9\textwidth]{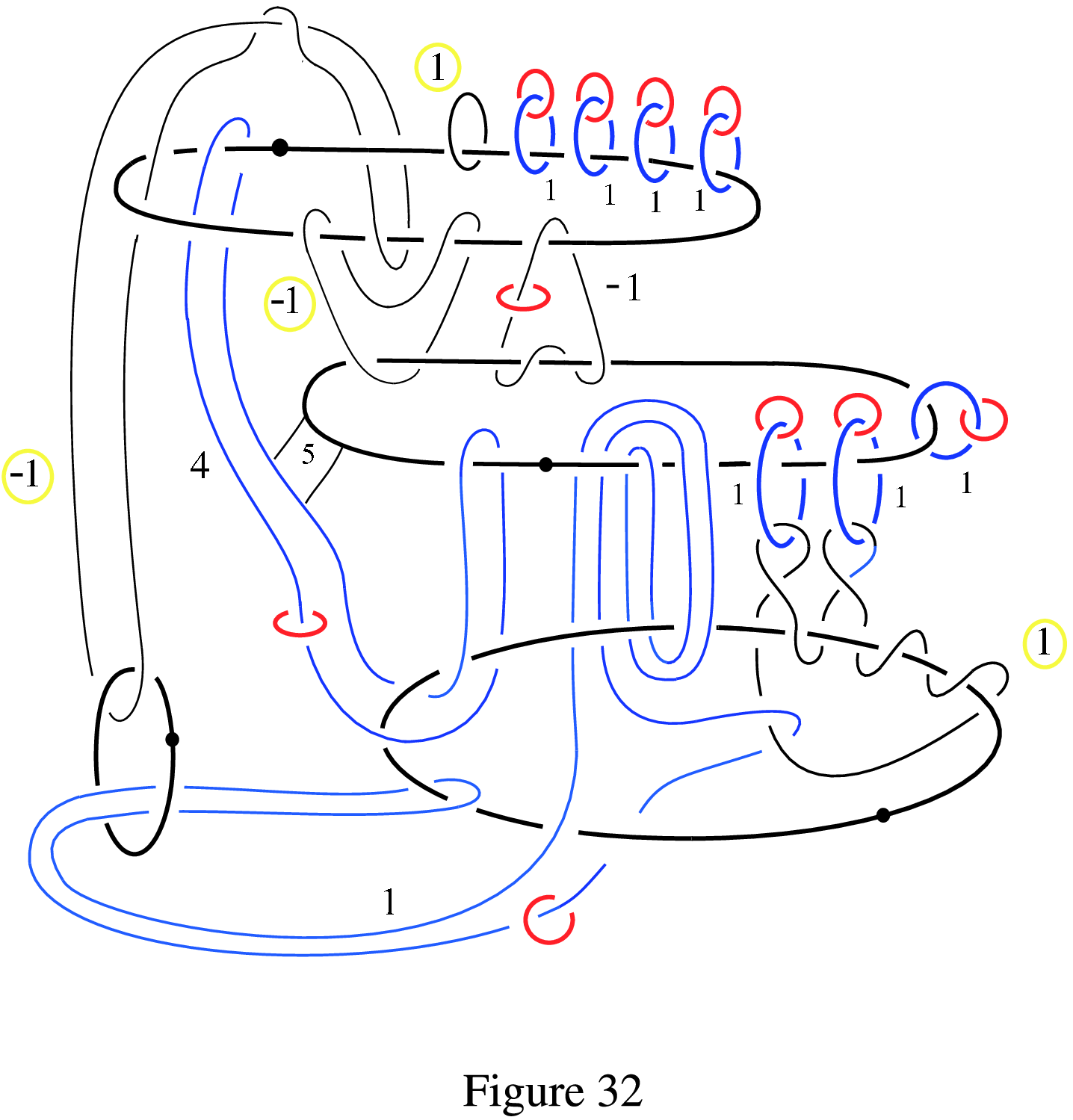}

\includegraphics[width=.9\textwidth]{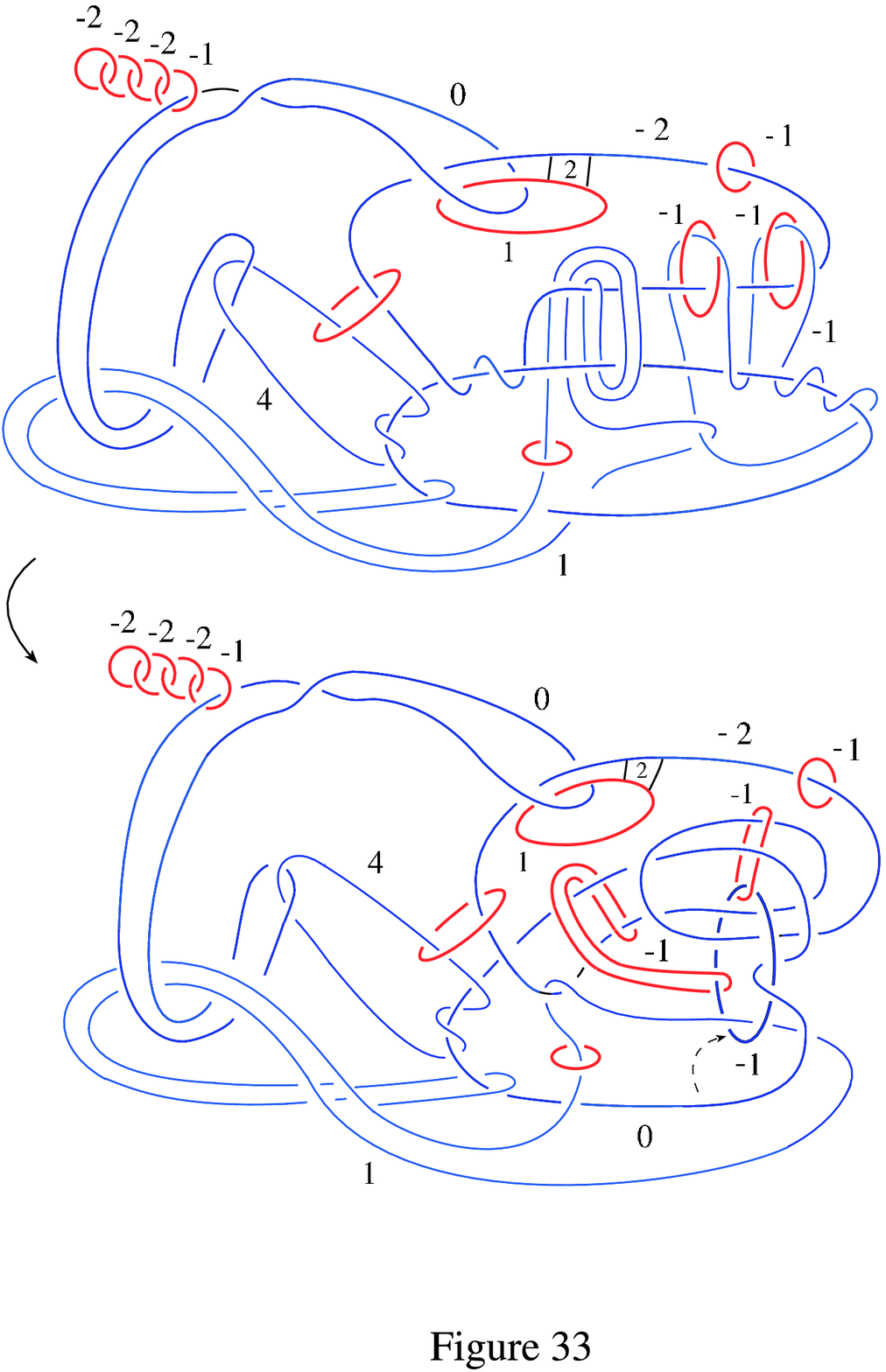}

\includegraphics[width=.8\textwidth]{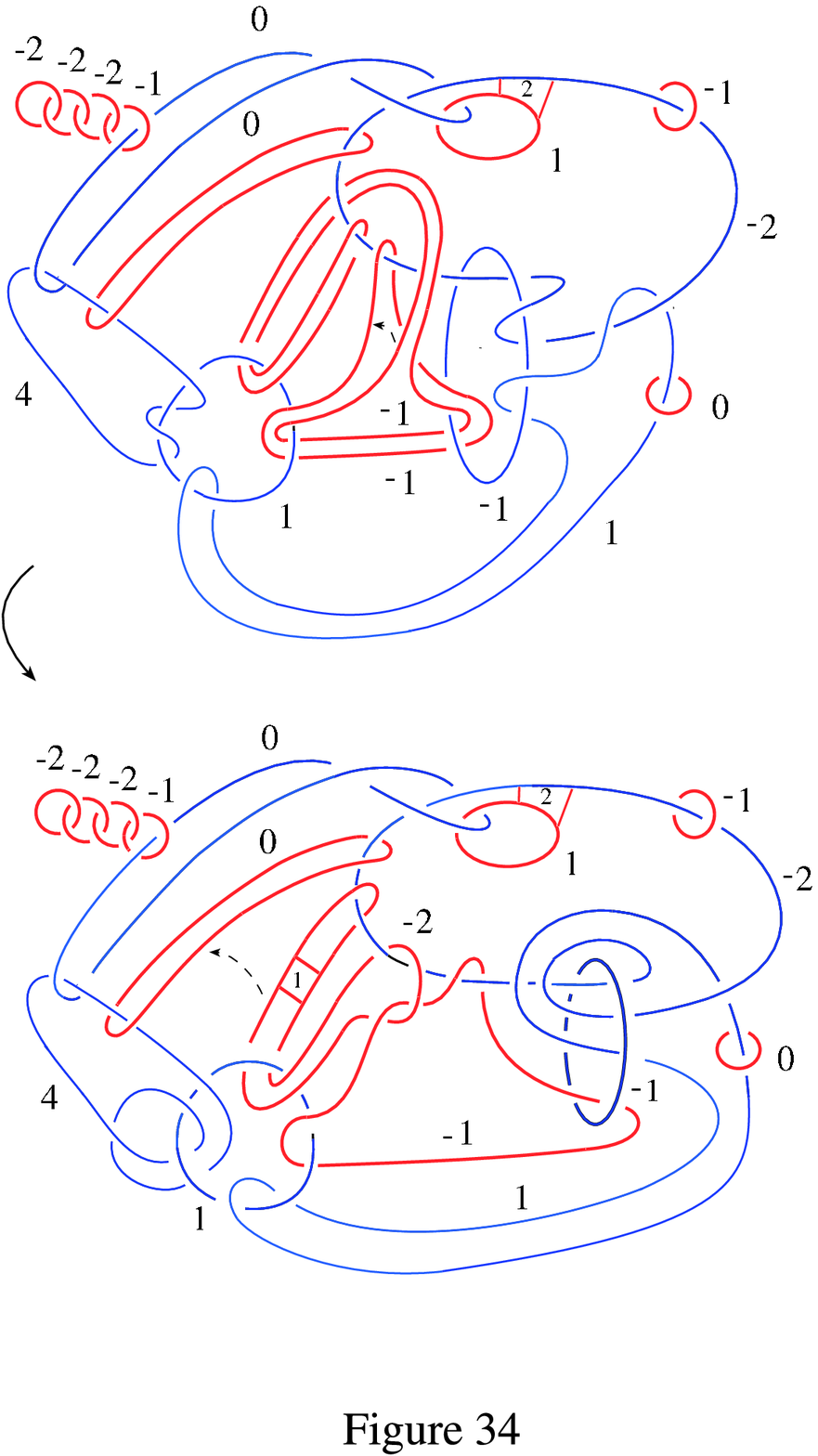}

\includegraphics[width=.85\textwidth]{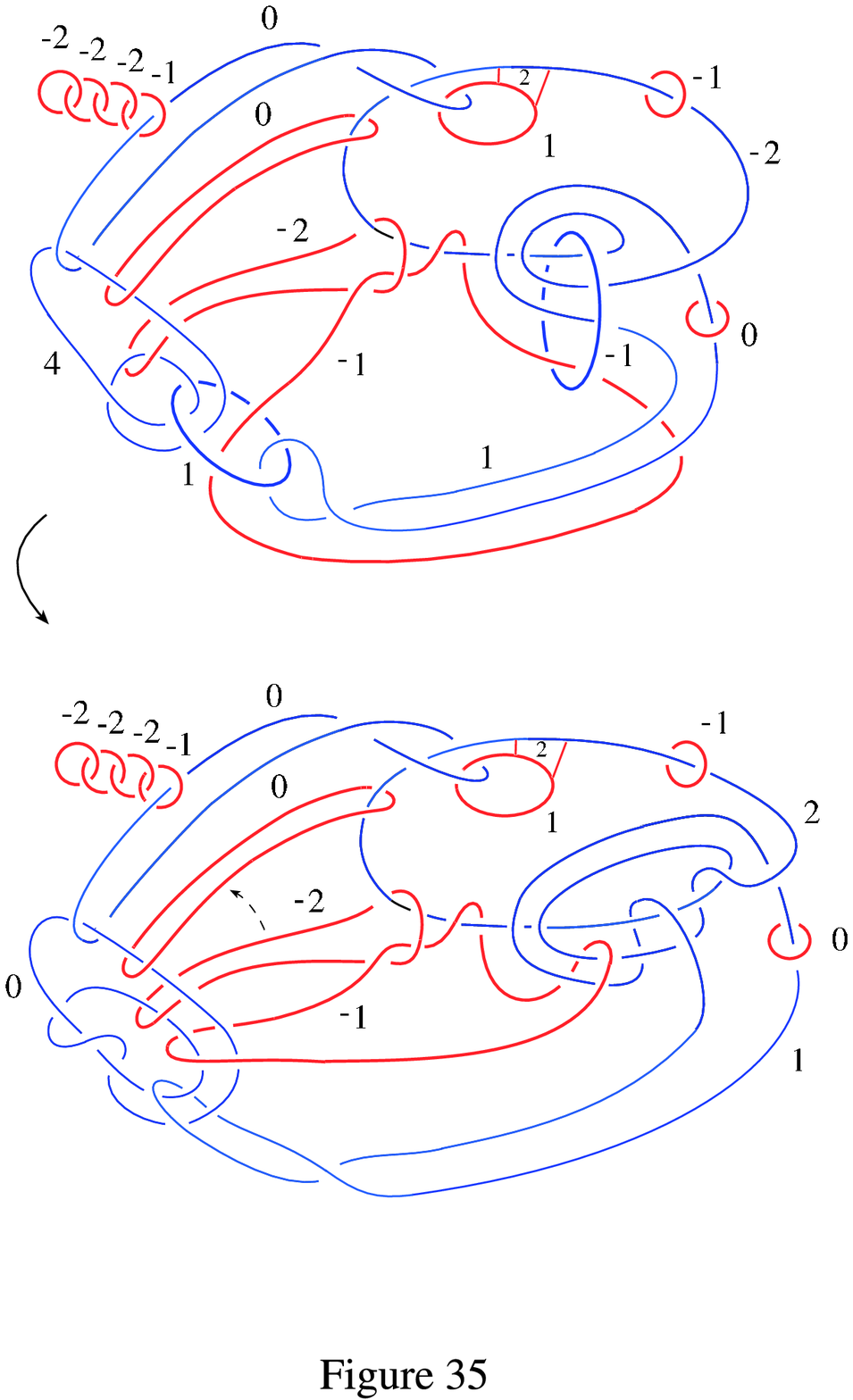}

\includegraphics[width=.9\textwidth]{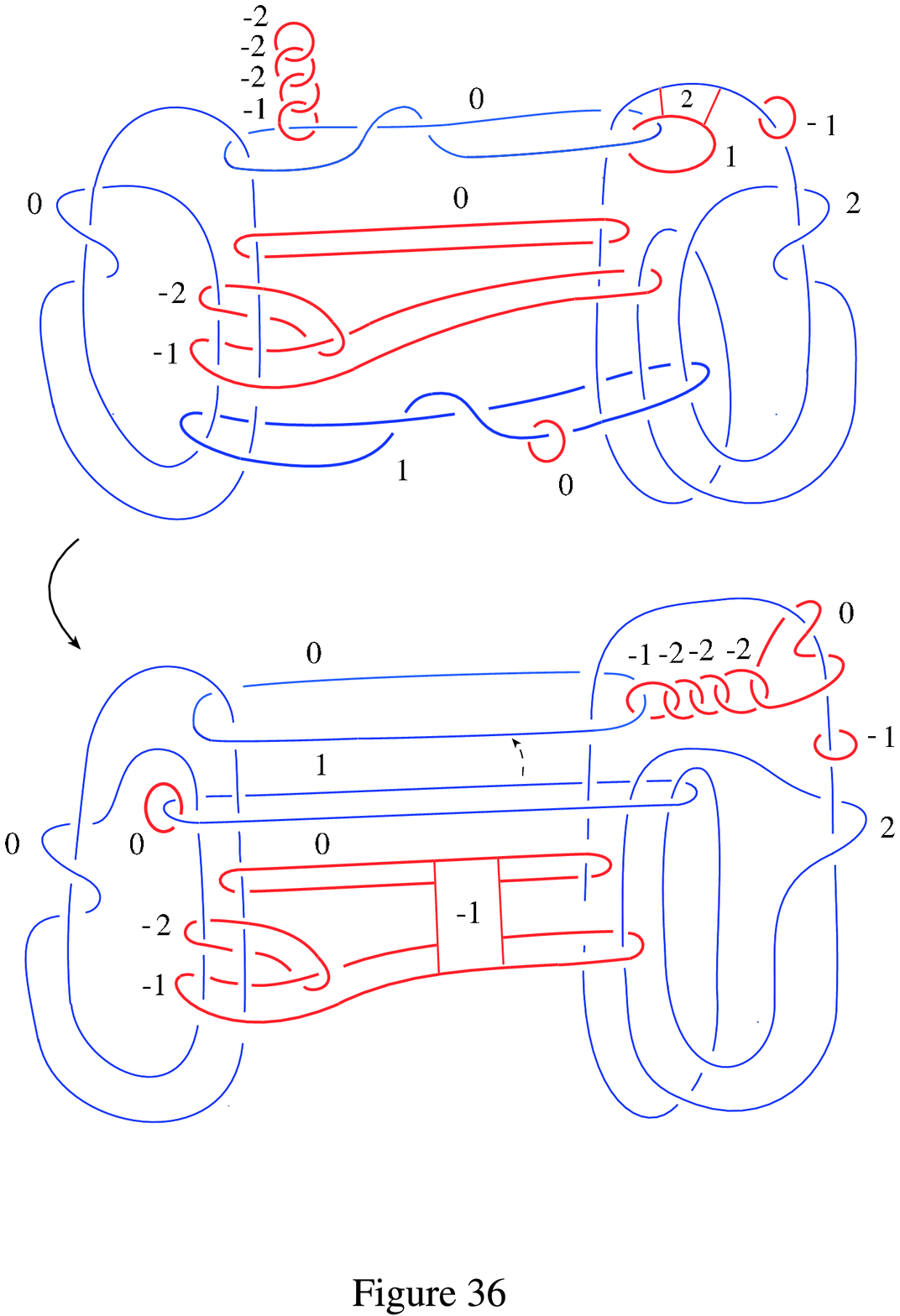}

\includegraphics[width=.9\textwidth]{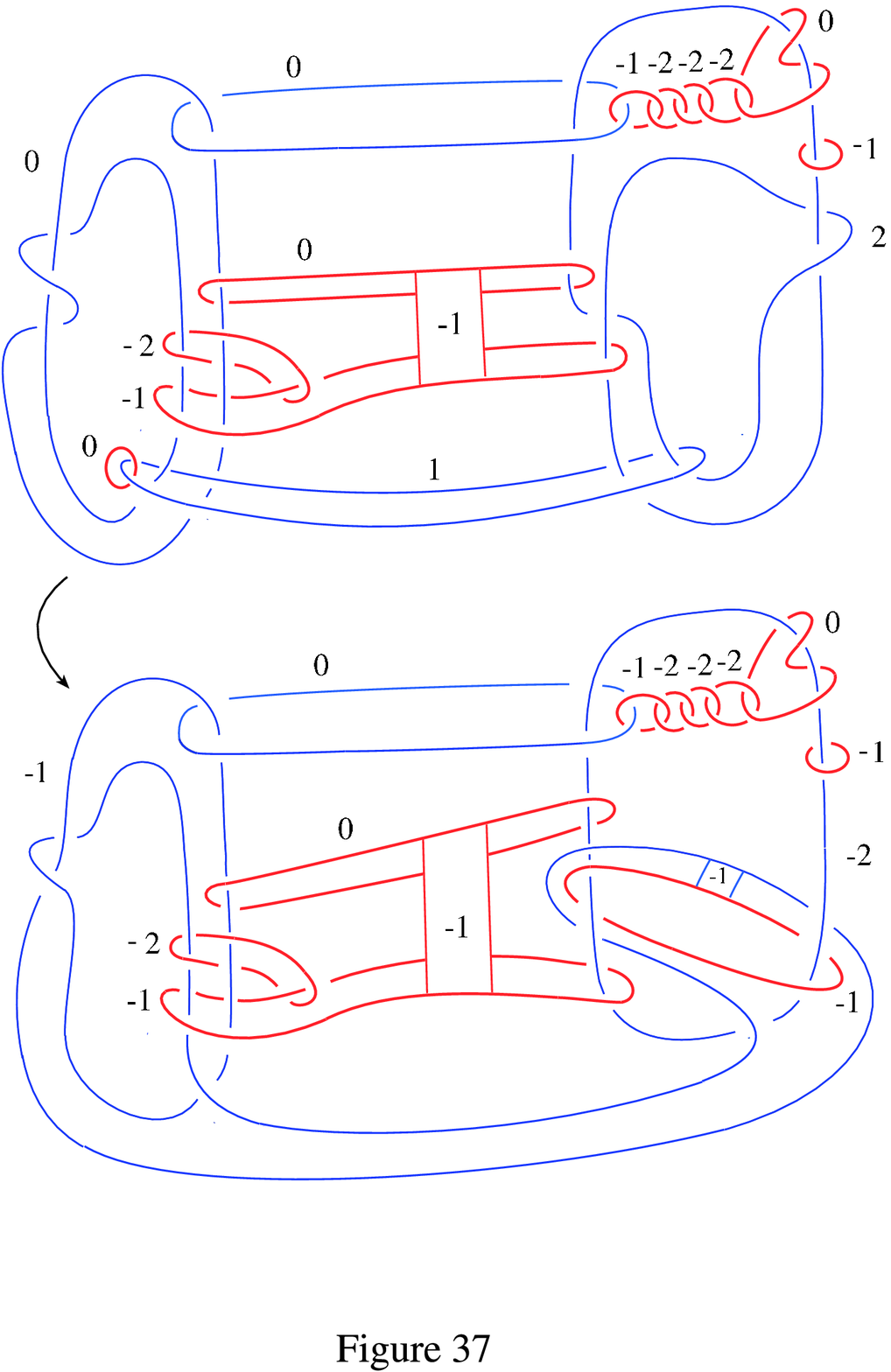}

\includegraphics[width=.9\textwidth]{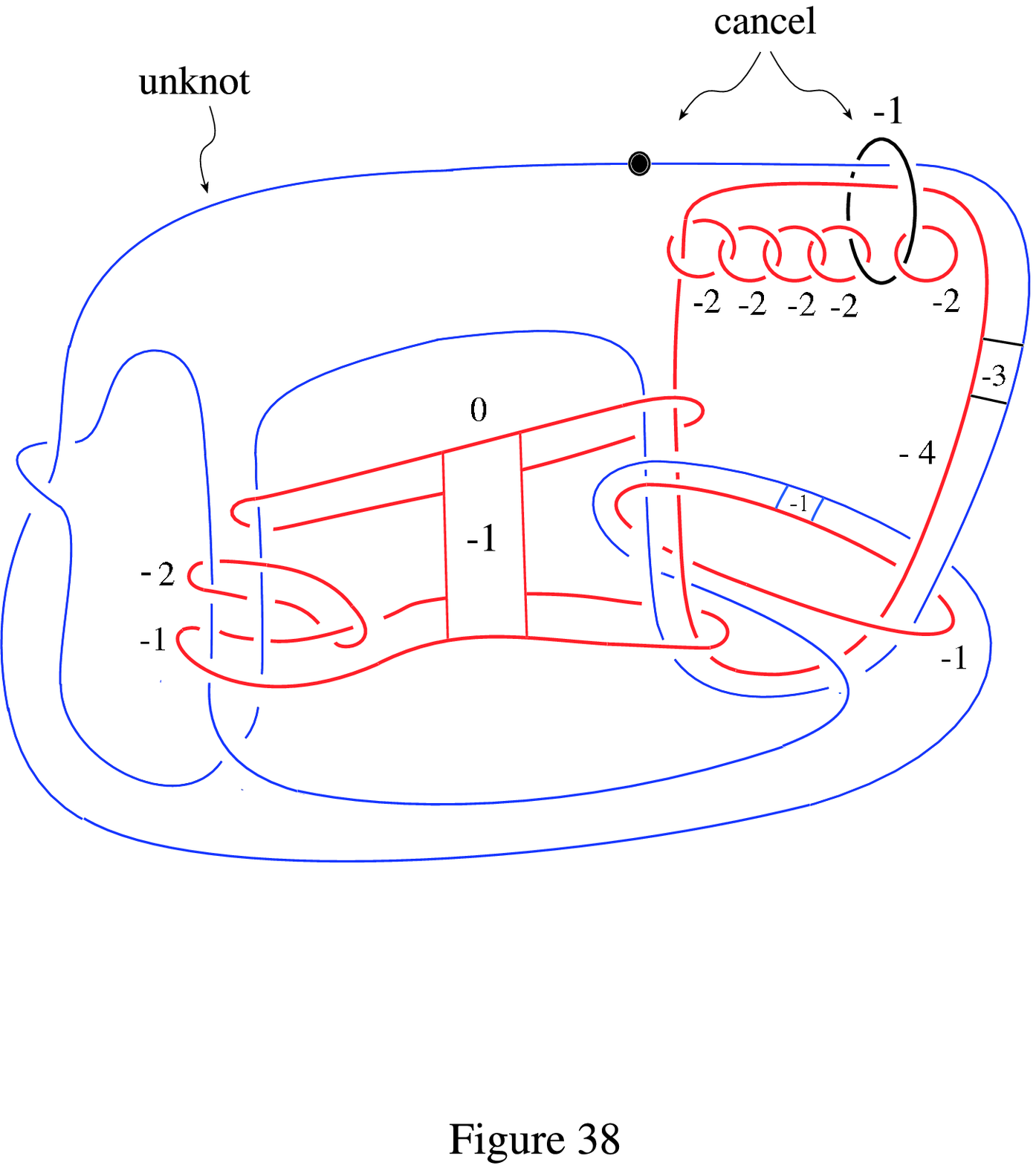}

\includegraphics[width=.75\textwidth]{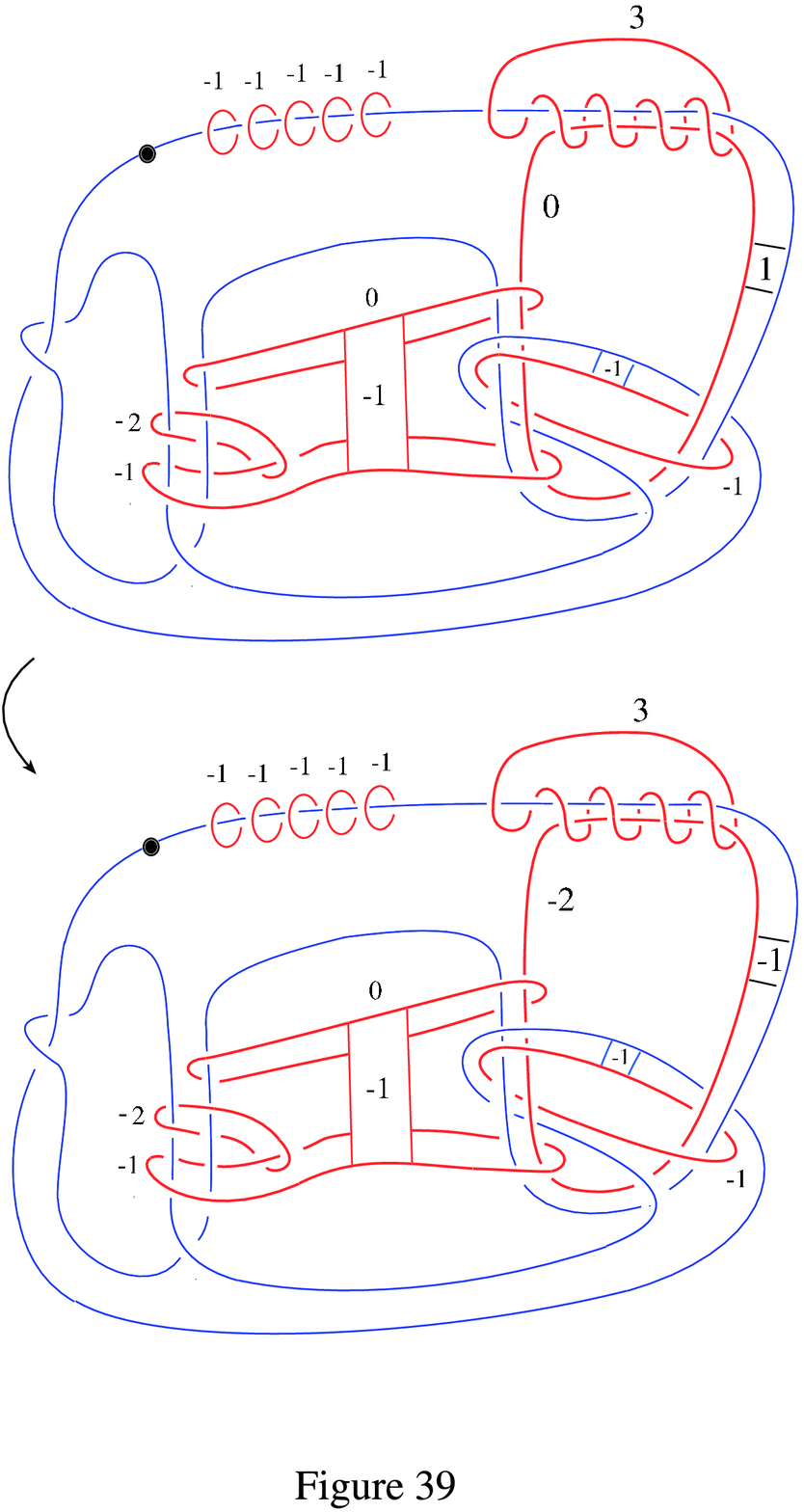}

\includegraphics[width=.75\textwidth]{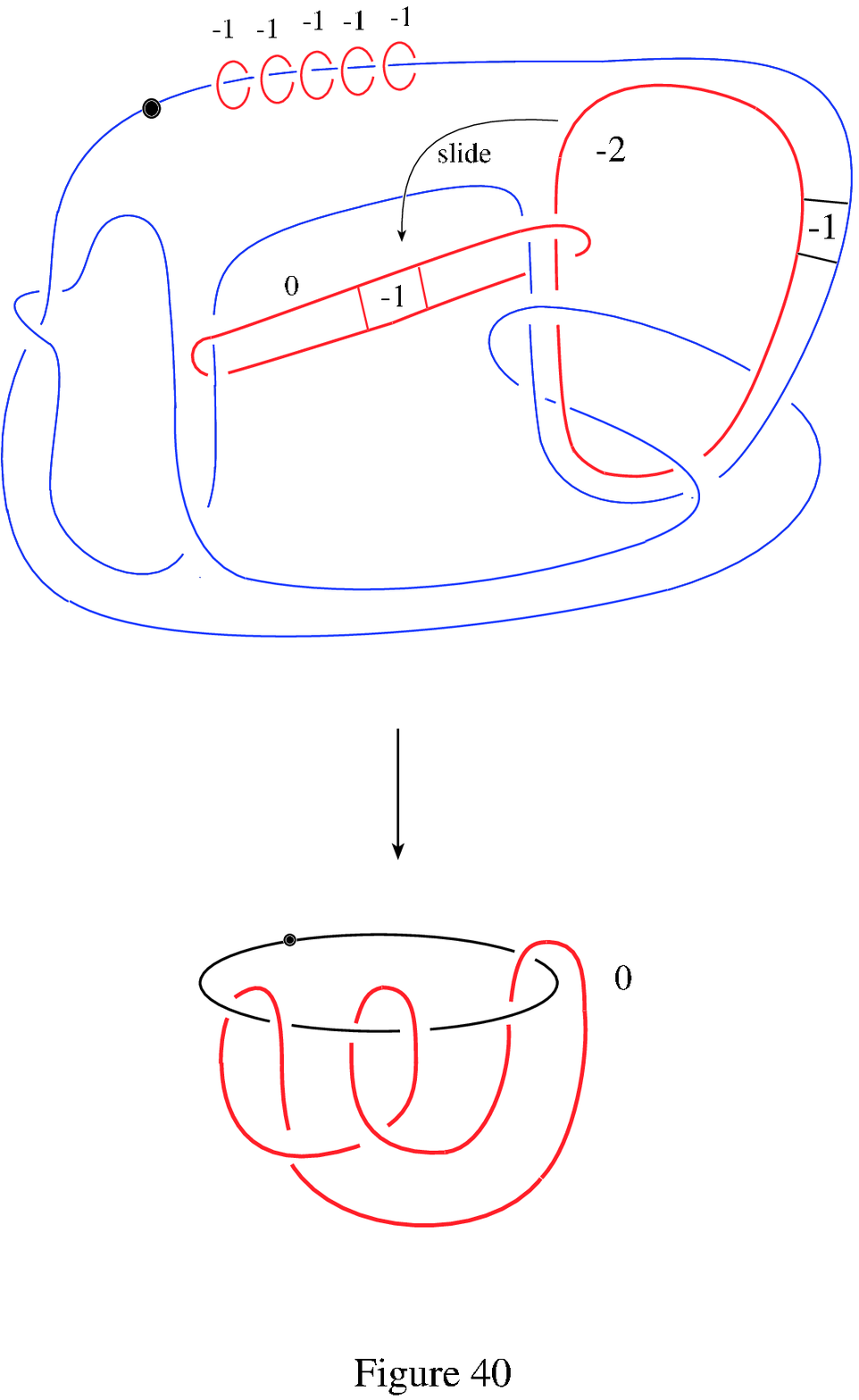}

\includegraphics[width=.75\textwidth]{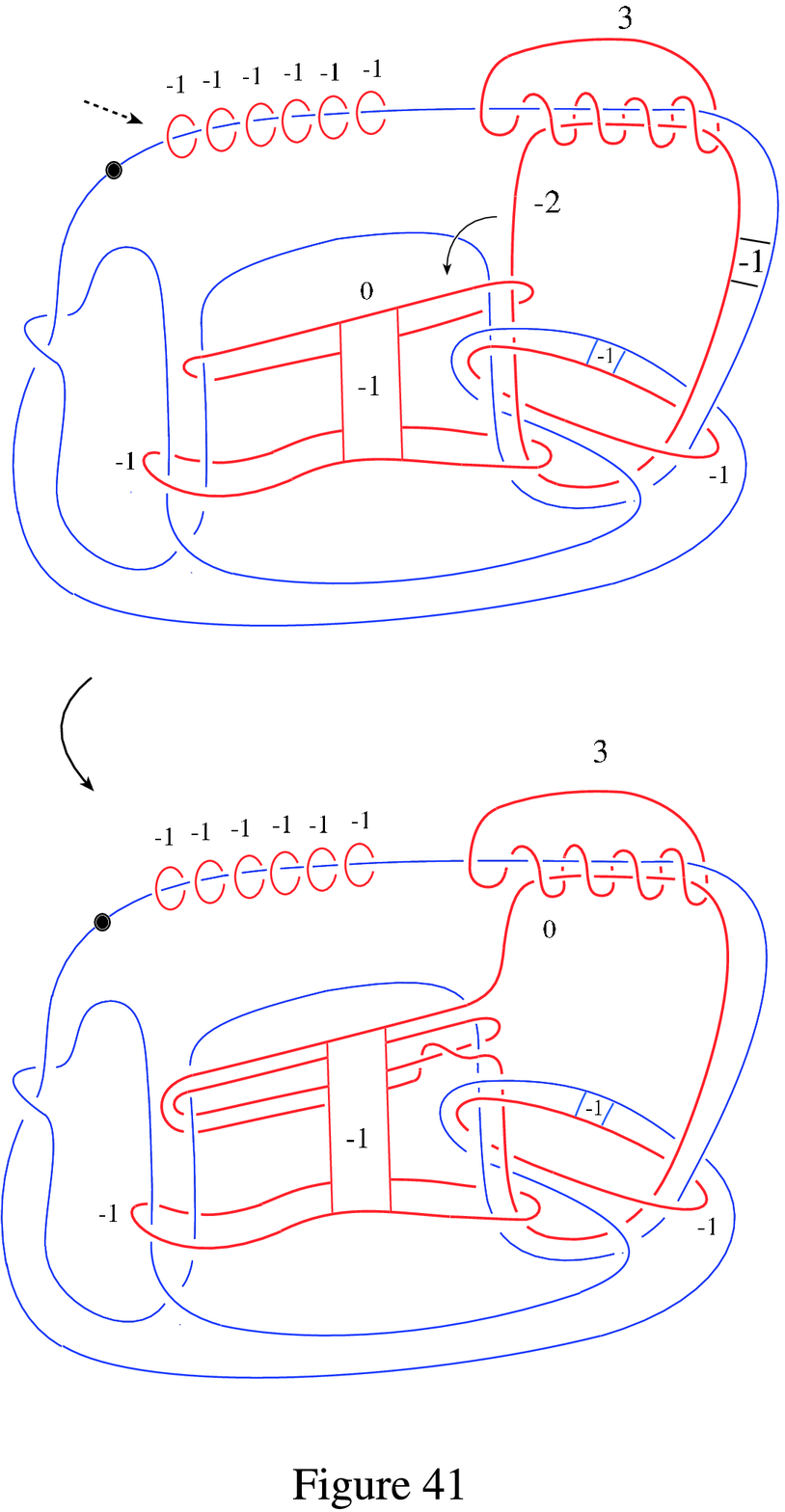}

\includegraphics[width=.9\textwidth]{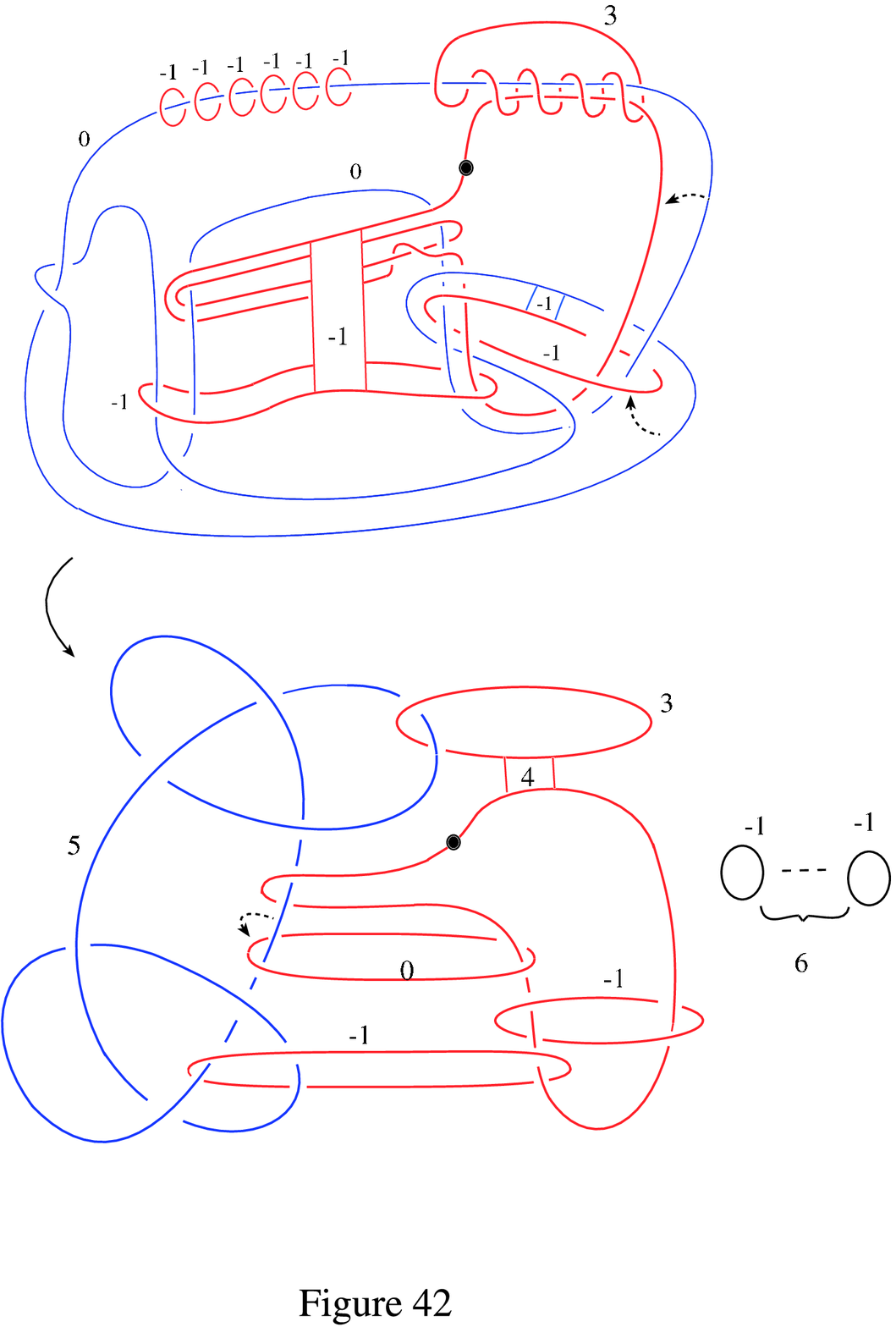}

\includegraphics[width=.9\textwidth]{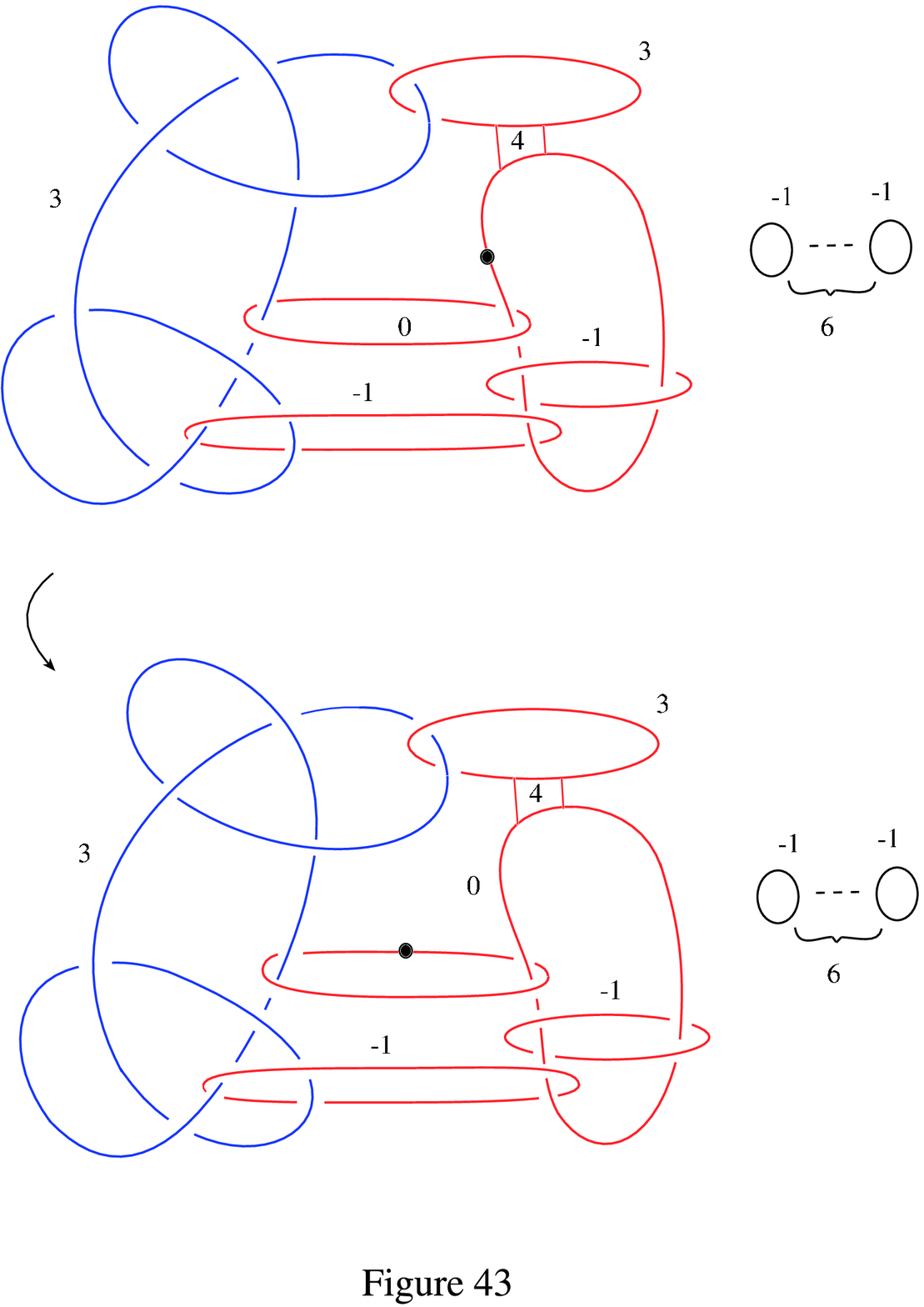}

\includegraphics[width=.9\textwidth]{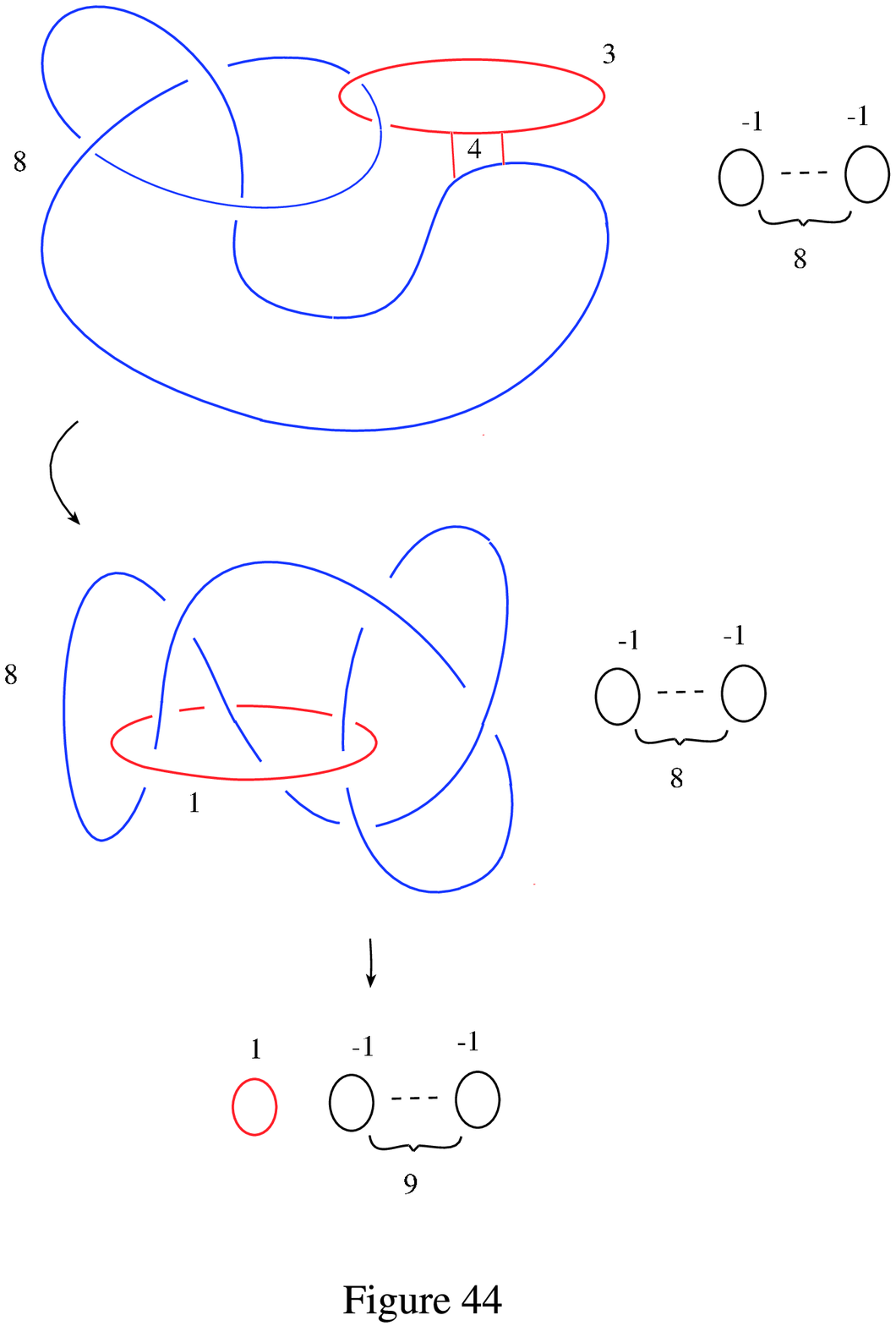}

\newpage

\begin{figure}[ht] 
   \begin{center}  
\includegraphics{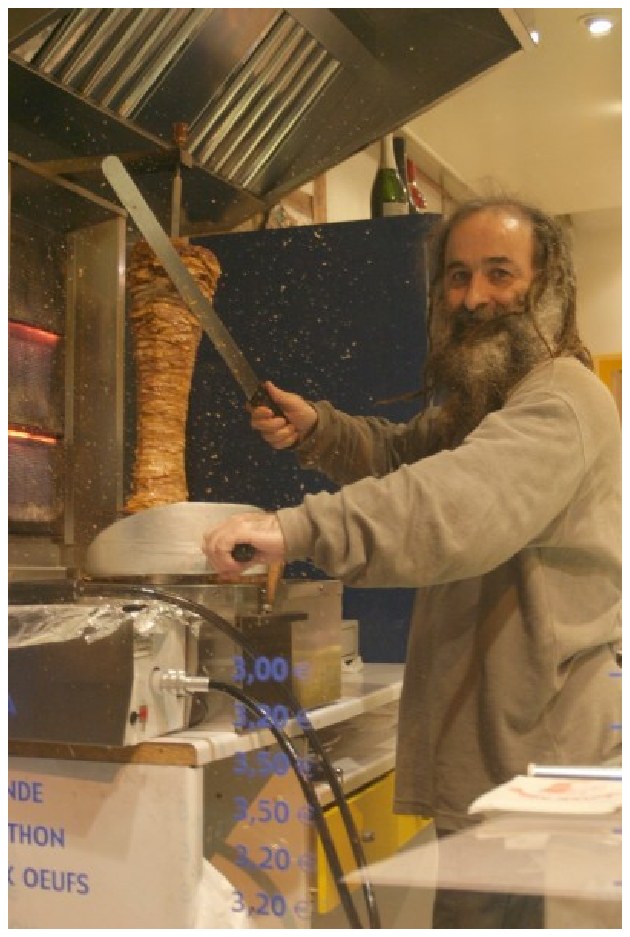}  
 \vspace{.1in}
 
 {\it If you keep turning handles of a 4-manifold upside down, while isotoping and canceling, you get a better picture of the manifold.}    
 \end{center}
\end{figure}

\end{document}